\documentstyle[12pt, my]{article}

\newfont{\Bbb}{msbm10 scaled \magstep1}

\newcommand\uc{\underline{\hbox{\Bbb C}}}

\newcommand\bC{\hbox{\Bbb C}}
\newcommand\bE{\hbox{\Bbb E}}
\newcommand\bR{\hbox{\Bbb R}}
\newcommand\bZ{\hbox{\Bbb Z}}
\newcommand\bX{\hbox{\Bbb X}}

\newcommand\bH{\hbox{\Bbb H}}

\newcommand\bS{\hbox{\Bbb S}}

\newfont{\es}{eusm10 scaled \magstep1}
\newfont{\ses}{eufm8 scaled \magstep1}
\newfont{\gt}{eufb10 scaled \magstep1}
\newfont{\sg}{eufb8 scaled \magstep1}
\newfont{\goth}{eufb10 scaled \magstep2}

\newcommand{\re}{\hbox{\gt Re}}
\newcommand{\im}{\hbox{\gt Im}}

\newcommand{\gA}{\hbox{\gt A}}

\newcommand{\gc}{\hbox{\gt c}}
\newcommand{\sgc}{\hbox{\sg c}}

\newcommand{\gF}{\hbox{\gt F}}

\newcommand{\gf}{\hbox{\gt f}}
\newcommand{\gG}{\hbox{\gt G}}

\newcommand{\gK}{\hbox{\gt K}}
\newcommand{\gL}{\hbox{\gt L}}

\newcommand{\gM}{\hbox{\gt M}}

\newcommand{\gq}{\hbox{\gt q}}
\newcommand{\goq}{\hbox{\gt q}}

\newcommand{\dir}{\hbox{\es D}}
\newcommand{\dol}{\hbox{\goth d}}
\newcommand{\can}{\hbox{\es K}}

\newcommand{\dps}{\dot{\psi}}
\newcommand{\dbps}{\dot{\bar{\psi}}}
\newcommand{\dtd}{\dot{\tau}_{\delta}}
\newcommand{\dad}{\dot{a}_\delta}

\def\ve{\varepsilon}
\def\vfi{\varphi}

\def\b{\cal B}
\renewcommand{\c}{{\cal C}}

\newcommand{\w}{{\cal W}}

\def\ra{\rightarrow}

\def\be{\begin{equation}}
\def\ee{\end{equation}}

\def\lan{\langle}
\def\ran{\rangle}
\def\na{\nabla}
\newcommand{\nab}{\mbox{$\boldmath{\nabla}$}}
\newcommand{\naf}{{}^\flat\nabla}
\newcommand{\onaf}{{}^\flat\overline{\nabla}}
\newcommand{\obnaf}{{}^\flat\overline{\mbox{$\boldmath{\nabla}$}}}

\def\nah{\hat{\nabla}}

\def\uu{\underline{u}}

\def\uso{\underline{so}}

\newtheorem{theorem}{Theorem}[section]
\newtheorem{lemma}[theorem]{Lemma}
\newtheorem{proposition}[theorem]{Proposition}
\newtheorem{corollary}[theorem]{Corollary}
\newtheorem{remark}[theorem]{Remark}
\newtheorem{definition}[theorem]{Definition}

\begin{document}

\title{Adiabatic limits of the Seiberg-Witten equations on  Seifert manifolds}

\author{Liviu I.Nicolaescu\thanks{{\bf Current address}: Dept.of Math.,University of Michigan, Ann Arbor, MI 48109-1003 USA ; liviu@math.lsa.umich.edu}}

\date{}

\maketitle

\addcontentsline{toc}{section}{Introduction}

\begin{center}
{\bf Introduction}
\end{center}

\bigskip

The  ``old''  instanton theory naturally lead to the instanton Floer homology of
 a 3-manifold $N$ as the missing piece in a general gluing formula for the     
 Donaldson invariants. Similarly, the Seiberg-Witten theory leads to the      
 ``monopole homology'' which is the Floer homology of the  Seiberg-Witten     
 functional  defined by a ${\rm spin}^c$ manifold (see [KM], [Mar], [MW], [Wa]).

The first main difficulty in understanding this homology comes from the fact that its defining 
chains   are not as explicit as the chains which
generate the instanton homology. In the latter case these are the flat connections
on an $SU(2)$ bundle  which are well understood both topologically and geometrically.   The
meaning of  the chains in monopole homology is far from obvious and the only explicit computations were made  when the 3-manifold $N$ is  a product $S^1 \times \Sigma$ where $\Sigma$ is a surface of genus $\geq 2$ (see [D] or [MST]).    The equations are  tractable 
 in this case is because  $S^1 \times N$ is a {\em K\"{a}hler manifold}. 
  As was pointed out in  [D], the solutions of the 3D  Seiberg-Witten equations coincide with the $S^1$-invariant solutions of  the 4D Seiberg Witten equations.  Fortunately, on a K\"{a}hler manifold  the solutions of these equations   can be described explicitly.

If  now  $N$ is  the total space of a principal $S^1$ bundle of {\em nonzero} 
degree  over a surface $\Sigma$ then $S^1 \times N$ admits a natural complex 
structure  but this time  the manifold  cannot be K\"{a}hler for the simple reason that the first Betti number is odd.

We   analyze the Seiberg-Witten equations on a special class of 
3-manifolds, namely those which admit a Killing vector field of constant  pointwise length and  satisfy an additional technical condition .   Topologically, these manifolds must be  Seifert fibered manifolds. 

On such manifolds the Dirac operators have an especially nice form and in particular  the Seiberg-Witten equations can be further dissected. We are interested in the behavior of the solutions of the Seiberg-Witten equations as the metric degenerates in the direction of the Killing vector field.  This corresponds to collapsing the fibers of the Seifert fibration.

The paper is divided  into four  parts.  The first part  studies in detail 
the differential geometry  of the metric almost contact (m.a.c.) manifolds.  
In particular, we distinguish a special class of such manifolds the so called {\em Killing m.a.c} manifolds.  These are Riemann manifolds which admit a Killing vector field $\zeta$ of pointwise length 1. Any oriented Killing m.a.c 3-manifold is diffeomorphic to a Seifert manifold and moreover any Seifert 3-manifold admits a Killing  m.a.c structure.   A special class of Killing  m.a.c. 3-manifolds consists of the  $(K, \lambda)$ manifolds, $\lambda\in {\bR}$.  They are characterized by the condition
\[
d\eta (x) \equiv 2\lambda \ast \eta
\]
where $\eta$ denotes the 1-form dual to the Killing vector field  $\zeta$ and $\lambda$ is a constant.  These manifolds are also characterized by the fact that their product with $S^1$ admits a natural integrable complex structure.

 The total  space of a principal $S^1$-bundle  admits a natural 
$(K, \lambda)$-structure   described for the first time  by Boothby and  Wang. 
 We present a 1-parameter family of  such nice  metric structures and we  
 explicitly compute its differential  geometric invariants: the Levi-Civita
  connection, the Ricci and the scalar curvature.   Factoring with  suitable 
   groups of isometries one can construct  many other interesting examples. 
   In particular, in  subsection $\S 1.4$ we show that any Seifert
    manifold  admits a natural $(K,\lambda)$-structure.  More precisely, the Thurston geometries  on  Seifert manifolds ( from the list of 6 described in [S])  are $(K,\lambda)$ structures. The scalar $\lambda$ is proportional with the Euler number of the Seifert fibration should be regarded as a measure of ``twisting '' of the fibration.  As the  metric degenerates (and so the fibers become shorter and shorter) $\lambda$ will  go to zero and thus the fibration will become ``less and less  twisted''.

As was observed by several authors ([ENS], [V]) the $(K,\lambda)$ structures with $\lambda >0$ are links of quasihomogeneous  singularities. Their 
Thurston geometry is  (almost) uniquely determined by  the analytical  structure of the singularity  and conversely, (see [Ne] or [Sch]) the Thurston geometry    fixes the analytical type of the
singularity.

The second part   is devoted to Dirac operators on m.a.c manifolds. 
 The spinor bundles corresponding to  the various  ${\rm spin}^c$ structures 
   can be very nicely described  in the almost-contact language. 
   We introduce two different classes of Dirac operators and compare them.  
   Also we establish a commutator identity which is an important ingredient  in the
  study of adiabatic limits.   In  the third part  we introduce the 3-dimensional  Seiberg-Witten
 equations and   study the adiabatic limits of solutions as the  metric 
 degenerates in the direction of $\zeta$.  On  an $S^1$- bundle of degree $\ell$ over  a surface of genus $g$ the ``adiabatic picture''   has some similarities  
   with the exact descriptions  in [D], [MST] or [Mun] when the 
 3-manifold is the trivial $S^1$-bundle over a surface.    However new phenomena occur when $|\ell|\geq 2g >0$. In this case,  for a certain range  of 
 ${\rm spin}^c$ structures  the adiabatic limit consists only of  the genuine  reducible solutions of the Seiberg-Witten equations which are adiabatic invariants.    We conclude this section by showing that if the fibers are sufficiently small  then  the reducible solutions are the only solutions of  the  Seiberg-Witten equations (corresponding to ${\rm spin}^c$-structures in the above range). They  determine a collection of   tori of  dimension $2g$ which  is nondegenerate in a Morse theoretic sense. 

The last section of the paper     illustrates how one can use the adiabatic knowledge  to actually produce  irreducible solutions of the Seiberg-Witten equations.  More precisely we show that the isolated, irreducible adiabatic solutions can always be perturbed to  genuine  solutions of the Seiberg-Witten equations corresponding to metrics with $\delta  \gg 1$.  The  idea is that these adiabatic solutions almost solve the Seiberg-Witten equations and if all goes well   the technique of  [T] can be used  to  detect nearby genuine irreducible solutions.   The problem  is technically   complex since all the important analytical quantities  such as eigenvalues, Sobolev constants etc. vary with  $\delta$    and one must   understand  quite accurately  the manner in which this happens.  There are fortunately many  magical  coincidences   which make this endeavour  less painful  than expected.

\tableofcontents

\section{Metric almost contact structures}

The tangent bundle of any oriented Riemann 3-manifold is trivial.  In particular,
 its structural group $SO(3)$ can  be reduced to $SO(2)\cong U(1)$. Such a reduction 
 is called an almost contact structure.  In this section we will discuss the special 
 differential geometric features of a 3-manifold  equipped with an almost contact 
 structure. Most of these facts are known (see [B], [YK])  but  we chose to present them in some detail in order to emphasize the special 3-dimensional features.

\subsection{Basic objects}
Consider an oriented  3-manifold $N$. 

\begin{definition}{\rm (a) An {\em almost contact structure} on $N$ is a  nowhere vanishing 1-form $\eta \in \Omega^1(N)$.

\noindent (b) A Riemann metric $g$ on $N$ is said to be compatible with  an almost 
contact structure  $\eta$ if  $|\eta (x)|_g=1$  for all $x\in N$. A {\em metric almost contact structure} 
(m.a.c) on $N$  is a pair $(\eta, g)$= (almost contact structure, compatible metric).}
\end{definition}

Consider a m.a.c  structure $(\eta, g)$ on the oriented 3-manifold $N$. 
 A local, oriented, orthonormal frame $\{\zeta_0, \zeta_1, \zeta_2\}$ of $TN$ 
 is said to be {\em adapted to the m.a.c. structure} if  $\zeta_0$ is the 
 metric  dual of  $\eta$.  The dual coframe of an adapted frame $\{\zeta_0, \zeta_1, \zeta_2\}$
   has the form $\{\eta^0, \eta^1, \eta^2\}$ where $\eta^0=\eta$ and $\ast \eta =\eta^2 \wedge \eta^2$. In the sequel we will operate exclusively with adapted frames.

Denote by $Cl(N)$ the bundle of  Clifford algebras generated by  $T^*N$ equipped with the induced metric.  The quantization map
\[
{\rm exterior\; \;algebra} \ra {\rm Clifford \;\; algebra}
\]
(see [BGV]) induces a map
\[
{\goq}:\Lambda^*T^*N \ra Cl(N).
\]
On the other hand $\Lambda^*T^*N$ has a natural structure of  $Cl(N)$-module so that  via the  quantization map we can construct an action of $\Lambda T^*M$ on itself
\[
{\gc}: \Lambda^*T^*N \ra {\rm End}\,(\Lambda^*T^*N)
\]
called {\em Clifford multiplication}.

On a m.a.c.  3-manifold  $(N,\eta, g)$  the Clifford multiplication by $\ast \eta$  has a remarkable property. More precisely 
\[
{\gc}(\ast \eta) \lan \eta\ran^\perp  =\lan\eta\ran^\perp \subset \Lambda^*T^*N.
\]   
If $(\eta^0, \eta^1, \eta^2)$ is a local coframe then the bundle $\lan \eta \ran^\perp$ is locally spanned by $\eta^1, \eta^2$ and ${\gc}(\ast \eta)$ acts according to  the prescription
\[
{\gc}(\ast \eta): \eta^1 \mapsto \eta^2\; ,\; \eta^2 \mapsto -\eta^1.
\]
In particular, we notice that both ${\gc}(\ast \eta)$ and $-{\gc}(\ast \eta)$ define complex structures on the real 2-plane bundle $\lan \eta\ran^\perp$.

\begin{definition}{\rm The complex line bundle $(\lan \eta \ran^\perp, -{\gc}(\ast \eta) )$ is 
called the {\em canonical line bundle} of the m.a.c. 3-manifold $(N,\eta, g)$ and is denoted by ${\can}={\can}_{\eta, g}$.}
\end{definition}

When viewed as a real bundle ${\can}$ (and hence ${\can}^{-1}$ as well) has a natural orientation.  We have an isomorphism of {\em oriented } real vector bundles
\be
T^*N \cong \lan -\eta \ran \oplus {\can} \cong \lan \eta \ran \oplus
{\can}^{-1}
\label{eq: candecomp}
\ee
where $\lan -\eta \ran$ (resp $\lan \eta \ran$)denotes the real line bundle 
spanned  and oriented by $-\eta$ (resp $\eta$).

\subsection{The structural equations of a m.a.c. manifold}

Consider an oriented  m.a.c. 3 manifold $(N, \eta ,g)$ and denote by  $\nabla$ the Levi-Civita connection of the metric $g$  Fix and adapted local frame of $\{\zeta_0, \zeta_1, \zeta_2\}$ and denote by $\{\eta^0, \eta^1, \eta^2\}$ the dual coframe. The connection 1-form of $\nabla$ with respect to these  trivializations can be computed using Cartan's structural equation.  More precisely if
\be
d \left[
\begin{array}{c}\eta^0 \\
\eta^1 \\
\eta^2
\end{array}
\right]
=\left[
\begin{array}{rrr}
0 & -A & B \\
A & 0 & -C \\
-B & C & 0 
\end{array}
\right] \wedge \left[
\begin{array}{c}\eta^0 \\
\eta^1 \\
\eta^2
\end{array}
\right]
\label{eq: structural}
\ee
($A$, $B$, $C$ are  real valued 1-forms locally defined on $N$) then
\be
\left\{
\begin{array}{lcrrr}
\nabla \zeta_0 & = &   & -A \otimes \zeta_1 & +B\otimes \zeta_2 \\

\nabla \zeta_1 & = & A \otimes \zeta_0& & -C\otimes \zeta_2 \\

\nabla \zeta_2 & = & -B \otimes \zeta_0& +C \otimes \zeta_1 & 
\end{array}
\right.
\label{eq: levicivita}
\ee
We will analyze the above equations when $(N, g)$ is a {\em Killing m.a.c.} manifold i.e. the vector $\zeta$ is Killing.  For $j=0,1,2$ we set
\[
A_j= i_{\zeta_j}A,\;\;B_j=i_{\zeta_j}B,\;\;C_j=i_{\zeta_j}C.
\]
Since $L_\zeta g=0$  we have the equality
\[
g(\nabla_X\zeta, Y)=-g(X,\nabla_Y\zeta)\;\;\;\forall X,Y\in {\rm Vect}\,(N).
\]
We substitute $X$ and $Y$ with pairs of basic vectors $\zeta_i$, $\zeta_j$  and we obtain the following  identities.

\be
A_0=B_0=A_1=B_2=0
\label{eq: kill1}
\ee
and 
\be
A_2(x)=B_1(x).
\label{eq: kill2}
\ee
Set  $\lambda (x):= A_2(x)=B_1(x)$. The structural equations  now yield
\[
d\eta = 2\lambda(x) \ast \eta.
\]
Thus the scalar $\lambda(x)$ is independent of the local frame used i.e. it is  an {\em invariant} of the  Killing m.a.c. structure  $(N, \eta, g)$.  

Differentiating the above equality we deduce  $0 =2 d(\lambda (x) \ast \eta)$ so that
\be
\partial_\zeta \lambda (x) =0.
\label{eq: constant}
\ee 
Note that
\[
\nabla_{\zeta}\zeta_1 =-C_0(x)\zeta_2,\;\;\;\nabla_\zeta\zeta_2 =C_0(x) \zeta_1.
\]
Thus $C_0(x)$ defines the infinitesimal rotation of $\lan \zeta \ran^\perp$ produced  by the parallel transport along $\zeta$.  Hence this is another invariant of the Killing m.a.c. structure and will be denoted by ${\vfi}(x)$.  Finally set
\[
b(x) =\lambda(x) + {\vfi}(x).
\]
Note for further references that
\be
[\zeta_1, \zeta_0]=\nabla_{\zeta_1}\zeta_0 -\nabla_0\zeta_1=b(x) \zeta_2
\label{eq: commu1}
\ee
and
\be
[\zeta_2, \zeta_0]= \nabla_{\zeta_2}\zeta_0 -\nabla_{\zeta_0}\zeta_2 =-b(x)\zeta_1.
\label{eq: commu2}
\ee

For each $\delta >0$ denote by $g_\delta$ the anisotropic deformation of $g$
defined by
\[
g_{\delta}(X,X)= g(X,X)\;\;{\rm if}\;\; g(X,\zeta) \equiv 0.
\]
\[
g_\delta (\zeta, \zeta)= \frac{1}{\delta^2}.
\]
Set 
\[
\eta_\delta = \eta/\delta.
\]
Note that  $d \eta_\delta =2\lambda \delta^{-1}\ast_\delta \eta_\delta$. Set $\lambda_\delta=\frac{\lambda}{\delta}$. Note that
\be
b_\delta = \delta b\;\;\;\; {\vfi}_\delta= \delta b-\lambda/\delta=\delta {\vfi}+\left(\delta -\frac{1}{\delta}\right)\lambda.
\label{eq: anisotrop}
\ee
Anisotropic deformations as above  were   discussed also in [YK] where they were named  D-homotheties.

\begin{remark}{\rm   (a) If $\lambda$ is not constant and $c$ is a regular value of $\lambda$ then the level set $\lambda^{-1}(c)$ is a smooth embedded surface in  $N$  and $\zeta $ is a nowhere  vanishing tangent vector field along $\lambda^{-1}(c)$. If $N$ is compact this  implies  $\lambda^{-1}(c)$ is an  embedded torus. The trajectories of $\zeta$ (which wander around the level sets $\lambda^{-1}(c)$) define an 1-dimensional Riemannian foliation on $N$.   Such  foliations on compact 3-manifolds where  completely classified in [Ca]. In particular,  the topological type of $N$ is severely restricted (see also Proposition \ref{prop: seifert} below).

\noindent (b) The above computations show that a Killing m.a.c structure defines a normal almost contact structure on $N$.  This means that the almost complex structure $J$ on $S^1\times N$ defined by
\[
J: \frac{\partial}{\partial \theta}\mapsto \zeta \;\;\;{\rm and} \;\;\;\zeta_1 \mapsto \zeta_2
\]
is integrable. See [B] or [YK] for more details.}
\end{remark}

We can now easily compute the sectional curvatures $ \lan R(\zeta_j, \zeta)\zeta, \zeta_j\ran $. More precisely we have
\[
R(\zeta_1,\zeta)\zeta=\left(\nabla_{\zeta_1}\nabla_\zeta -\nabla_\zeta\nabla_{\zeta_1}-\nabla_{[\zeta_1,\zeta]}\right)\zeta
\]
\[
=-\nabla_\zeta(\lambda(x)\zeta_2)-b(x)\nabla_{\zeta_2}\zeta=\lambda^2(x)\zeta_1.
\]
Hence
\be
\lan R(\zeta_1, \zeta)\zeta, \zeta_1\ran=\lambda^2(x)
\label{eq: curv1}
\ee
and similarly
\be
\lan R(\zeta_2, \zeta)\zeta, \zeta_2\ran  =\lambda^2(x).
\label{eq: curvi3}
\ee
Hence the scalar curvature of $N$ is  determined by 
\be
s=2\kappa +4\lambda^2(x)
\label{eq: scacurv}
\ee
where $\kappa(x)$ denotes the sectional curvature of the plane spanned  by $\zeta_1$ and $\zeta_2$. 

 The connection $\nabla$ defines via orthogonal projections a connection $\nabla^\perp$ on
 the complex line bundle ${\rm Ann}\,\eta =\lan \zeta\ran^\perp$. (The complex
 structure is defined by ${\bf i}\zeta_1=\zeta_2$.  More precisely
\[
\left\{
\begin{array}{lcrr}
\nabla^\perp \zeta_1& =&  &-C(x)\otimes \zeta_2\\
\nabla^\perp\zeta_2 &= & C(x)\otimes \zeta_1 & \\
\end{array}
\right.
\]
where $C(x)= {\vfi}(x)\eta +C_1(x)\eta_1+C_2(x)\eta_2$. Using the complex
structure in $\lan\zeta\ran^\perp$ we can locally describe $\nabla^\perp$ as
\[
\nabla^\perp =d-{\bf i} C(x).
\]
Under anisotropic deformations  this connection 1-form changes to
\[
C_\delta  ={\vfi}_\delta\eta_\delta +C_1\eta_2+C_2\eta_2 =\frac{1}{\delta}{\vfi}_\delta \eta +C_1\eta_2+C_2\eta_2.
\]
The equality (\ref{eq: anisotrop}) shows that as $\delta \ra \infty$ the form  $C_\delta$ converges to
\[
C_\infty:=b(x)\eta +C_1(x)\eta_1+C_2(x)\eta_2 =\lambda(x)\eta +C(x).
\]
If we denote by $\nabla^{\infty}$ the limiting connection we have
\be
\nabla^\infty =\nabla^\perp -{\bf i}\lambda(x)\eta.
\label{eq: anisotrop1}
\ee
Denote by $F^\perp$ the curvature of the connection  $\nabla^\perp$ and by $\sigma$ the scalar
\be
\sigma=\lan F^\perp(\zeta_1,\zeta_2)\zeta_2, \zeta_1 \ran.
\label{eq: defcurv}
\ee
It is not difficult to show that $\sigma$ is independent of the local frame and so is an invariant of the  Killing m.a.c. structure.    With respect to this frame it has  the description
\[
\sigma(x)= \partial_{\zeta_1}C_2 -\partial_{\zeta_2}C_1 -(C_1)^2-(C_2)^2 .
\]
Note that the anisotropic deformation introduced in $\S 1.1$  does not change $\nabla^\perp$ so that
\[
\sigma_\delta(x)=\sigma (x).
\]
Using the structural equations of $\nabla$ we deduce
\[
\nabla_{\zeta_1}\nabla_{\zeta_2}\zeta_2 =(\partial_{\zeta_1}C_2)\zeta_1 -C_1C_2 \zeta_2.
\]
\[
\nabla_{\zeta_2}\nabla_{\zeta_1}\zeta_2 =(\lambda C_1-\partial_{\zeta_2}\lambda)\zeta +\left(\lambda^2(x) \zeta_1 +\partial_{\zeta_2}C_1\right)\zeta_1  -C_1C_2\zeta_2
\]
\[
[\zeta_1,\zeta_2]=-2\lambda(x)\zeta +C_1\zeta_1 +C_2\zeta_2.
\]
Hence
\[
\kappa(x) =\lan R(\zeta_1,\zeta_2)\zeta_2, \zeta_1\ran 
\]
\[
=\partial_{\zeta_1}C_2 -\partial_{\zeta_2}C_1 -(C_1)^2-(C_2)^2 -\lambda^2(x) +2\lambda(x){\vfi}(x)
\]
\be
=\sigma(x) -\lambda^2(x) +2 \lambda(x) {\vfi}(x).
\label{eq: curv4}
\ee
In particular, using (\ref{eq: scacurv}) we deduce
\[
s(x) = 2\{\sigma(x) +\lambda^2(x)  +  2\lambda (x){\vfi}(x)\}.
\]
Note that
\[
s_\delta(x)= 2\{\sigma(x) +\lambda^2/\delta^2 +2\lambda(b-\lambda/\delta)/\delta\}
\]
so that
\be
\lim_{\delta \ra \infty} s_\delta(x) = 2\sigma(x).
\label{eq: limit}
\ee
Thus for $\delta$ very large $\sigma(x)$ is a good approximation for the scalar curvature.

\begin{proposition}{\rm Assume $N$ is a compact oriented  Killing m.a.c 3-manifold.  Then $N$ is diffeomorphic to an oriented Seifert fibered 3-manifold. Conversely, any compact oriented Seifert 3-manifold admits a Killing m.a.c structure.}
\label{prop: seifert}
\end{proposition}

\noindent{\bf Proof} \hspace{.3cm} The Killing  m.a.c. structures $(\eta, g)$  with respect to the {\em fixed} metric $g$ are parameterized by the unit sphere in the Lie algebra of compact Lie group of isometries ${\rm Isom}\,(N,g)$.  In particular, the group ${\rm Isom}\,(N, g)$ has positive dimension. If this is the case the maximal torus containing $\zeta$ induces at least one fixed-point-free $S^1$ action  on $N$  (slight perturbations of $\zeta$ in the Lie algebra of ${\rm Isom}(N,g)$ will not introduce  zeroes of the corresponding vector field on $N$).  Hence $N$ must be a Seifert manifold. 

Conversely, given a Seifert fibered  manifold $N$  denote  by $\zeta$ the generator of the fixed-point-free $S^1$ action and for each $\theta \in S^1$ denote by ${\cal R}_\theta$ its action on $N$   Define ${\gM}_\zeta$ as the collection of Riemann metrics $g$ on $N$ such that $|\zeta(x)|_g\equiv 1$. Note that   ${\gM}$ is convex and
\[
{\cal R}_\theta^*{\gM}_\zeta \subset {\gM}_\zeta.
\]
The $S^1$-average  of any $g\in {\gM}_\zeta$  is $\zeta$ -invariant and thus defines a Killing m.a.c structure on $N$. \hspace{.3cm} $\Box$

\bigskip

When the invariant $\lambda(x)$ is constant  we will call the structure  $(N,\eta, g)$ a  $(K, \lambda)$-manifold.  The $(K,1)$-manifolds are also known as $K$-{\em contact} manifolds (cf. [B], [YK]). In dimension 3 this notion coincides with the notion of Sasakian
manifold.

On a $(K, \lambda)$ manifold  the sectional curvature of  any plane containing $\zeta$ 
 is $\lambda^2$  and the full curvature
  tensor is completely determined by the the sectional curvature $\kappa$ of the planes 
  spanned by $\zeta_1$ and $\zeta_2$.   The scalar curvature of  $N$ is
\[
s= 2\kappa +4\lambda^2.
\]
 Using the structural equations we deduce  that  with respect to the adapted frame $\{\zeta_0, \zeta_1, \zeta_2\}$ the Ricci 
curvature has the form
\[
{\rm Ric}=\left[
\begin{array}{rrr}
2 \lambda^2 &  0& 0 \\
0 &\kappa +\lambda^2 & -\kappa\\

0 & -\kappa &\kappa +\lambda^2 
\end{array}
\right]
\]

For any oriented $(K,\lambda)$ 3-manifold  (not necessarily compact) we denote by ${\gK}_N$ the group of isomorphism of the $(K,\lambda)$ structure i.e. orientation preserving isometries which invariate $\eta$.  For any discrete subgroup   $\Gamma \subset {\gK}_N$ acting freely and discontinuously  on $N$ we obtain a covering 
\[
N \ra N/\Gamma.
\]
Clearly $N/\Gamma$ admits  a natural $(K, \lambda)$ structure  induced from  $N$.   In $\S 1.4$ we will use this  simple observation to construct $(K, \lambda)$-structures on any Seifert 3-manifold.

\subsection{The Boothby-Wang construction}

 In this subsection we describe some natural  $(K,\lambda)$ structures on the
total space of  a principal $S^1$-bundle over a compact oriented surface.
Except some minor modifications this construction is due to Boothby and  Wang,
[BW] (see also [B]).

Consider $\ell \in {\bZ}$ and  denote by $N_\ell$ the total space of a degree
$\ell$ principal $S^1$ bundle over a  compact oriented surface of genus $g$: $S^1 \hookrightarrow N_\ell \stackrel{\pi}{\ra} \Sigma.$ We orient $N_\ell$ using the rule
\[
\det T N_\ell = \det TS^1 \wedge \det T \Sigma.
\]
Assume $\Sigma$ is equipped with a Riemann metric $h_b$ such that ${\rm vol}_{h_b}(\Sigma)= \pi$. Recall that if $\omega \in {\bf i}{\bR}\otimes \Omega^1(N_\ell)$ is a
connection form on $N_\ell$ then $d\omega$ descends
to a 2-form $\Omega$ on $\Sigma$: the curvature of $\omega$. Moreover
\[
\int_\Sigma c_1(\Omega)=\frac{\bf i}{2\pi} \int_\Sigma \Omega = \ell = \int_\Sigma
\frac{\ell}{\pi }dv_{h_b}.
\]
Notice in particular that ${\bf i}\Omega$ is cohomologous to  $2\ell dv_{h_b}$.
 Now pick a connection form $\omega$ such that
 \[
 {\bf i }\Omega = 2 \ell dv_{h_b}.
 \]
 Denote by $\zeta$ the unique vertical vector field on $N_\ell$ such that
 $\omega (\zeta)\cong {\bf i}$ and set $\eta_\delta = -{\bf i} \omega/\delta$.   Thus
 $\eta_\delta (\zeta) \cong 1/\delta$. Notice that ${\rm Ann}\, \eta$ coincides with the horizontal
 distribution $H_\omega$ defined by the connection $\omega$.   Now define  a
 metric $h$ on  $N_\ell$  according to the prescriptions.
  
 \[
 h(\zeta, \zeta)= 1/\delta^2
 \]
 and
 \[
 h(X, Y)= h_b(\pi_*X, \pi_*Y)\;\;{\rm if}\; X,\; Y\; {\rm are\; horizontal}.
 \]
 Clearly $L_\zeta h=$ i.e. $\zeta $ is a Killing vector field.  Moreover
 \[
 \delta d\eta_\delta  =- \Omega =- 2 \ell \pi^* dv_{h_b} = -2\ell \ast
 \eta_\delta.
 \]
 In other words, we have constructed a $(K, - \ell/\delta)$ structure on $N_\ell$.
 
 We conclude this subsection by describing  the geometry of $N_\ell$ in terms
 of the geometry of $\Sigma$.  The only thing we need to determine is the
 curvature $\kappa$ introduced in the previous subsection in terms of the
 sectional curvature $\Sigma$.   It is not difficult to see this coincides with the  invariant $\sigma(x)$ defined in (\ref{eq: defcurv}). Because of the special geometry of this situation formula (\ref{eq: curv4}) can be  further simplified.

 To achieve this we will use again the structural equations. Denote by
 $\{ \psi^1, \psi^2\}$ a local, oriented orthonormal coframe on $\Sigma$ and
 set $\eta^j =\pi^* \psi^j$, $j=1,2$. Then $\{\eta^0=\eta, \eta^1, \eta^2\}$ is
 an adapted  coframe on $N_\ell$. 
 
 The structural equations of $\Sigma$ have the form
 \be
 d \left[
 \begin{array}{r}
 \psi^1 \\
 \psi^2
 \end{array}
 \right]=
 \left[
 \begin{array}{rr}
 0 & -\theta \\
 \theta & 0 
 \end{array}
 \right] \wedge \left[
 \begin{array}{r}
 \psi^1 \\
 \psi^2
 \end{array}
 \right]
 \label{eq: structuralb}
 \ee
where $\theta$ is a 1-form locally defined on $\Sigma$.  Set
\[
\theta_j =i_{\pi_*\zeta_j}\theta,\;\;j=1,2
\]
Since $ d\eta^j = \pi^* d \psi^j$ we deduce that $d\eta^j$ is horizontal. On
 the other, hand using  (\ref{eq: structural}) we deduce
 \[
 d\eta^1 =A\wedge \eta^0 -C\wedge \eta^2  =(A_2 + C_0)\eta^2 \wedge \eta^0 -C_1 \eta^1 \wedge \eta^2.
 \]
 This implies  $C_0 = -A_2 = \ell/\delta$ so that
\be
b(x)\equiv 0
\label{eq: commute}
\ee
and
\be {\vfi} \equiv   \ell /\delta.
 \label{eq: bw1}
 \ee
 Using  (\ref{eq: structuralb}) we deduce $d\psi^1 = -\theta \wedge  \psi^2 = -\theta_1 \psi^1 \wedge \psi^2$ which  yields
 \be
 C_1= \theta_1.
\label{eq: bw2}
\ee
Similarly one shows 
\be
C_2=\theta_2.
\label{eq: bw3}
\ee
Using the above two equalities  and (\ref{eq: anisotrop1}) we conclude that the limiting connection  $\nabla^{\infty}$ on ${\can}$ is none-other than the pullback by $\pi:N\ra \Sigma$ of the Levi-Civita connection on the canonical line bundle $K_\Sigma\ra \Sigma$ i.e.
\be
\nabla^\infty =\pi^*\nabla^K.
\label{eq: anisotrop2}
\ee 
Using (\ref{eq: curv4}) and (\ref{eq: bw1}) we deduce
\be
\kappa =\sigma -3\ell^2/\delta^2
\label{eq: bw4}
\ee
and in particular the scalar curvature of $N_\ell$ is given by
\be
s_N = 2(\sigma - \ell^2/\delta^2).
\label{eq: bw5}
\ee

\begin{remark}{\rm  Let $N$  denote   the total space of a principal $S^1$ bundle over  an oriented surface $\Sigma$, {\em not necessarily compact}. If $\omega$ denotes a connection  form on $N$ such that   $-{\bf i} d\omega$ descends to a {\em constant} multiple of the volume form on  $\Sigma$ then  the previous computations extend verbatim to this case and one sees that in this situation one  also obtains  a $(K,\lambda)$-structure on $N$. 

For example let $N$ denote the unit tangent bundle of the hyperbolic plane  ${\bH}^2$. The Levi-Civita connection on ${\bH}^2$ induces an $S^1$-connection $\omega$ on $N$.   Then
\[
-{\bf i} d\omega = -1 d\,vol_{{\bH}^2}
\]
since ${\bH}^2$  has constant  curvature $\equiv -1$.  Thus $N$ has a natural $(K, 1)$ (= Sasakian) structure.

  The group of isometries of ${\bH}^2$ is  $PSL(2, {\bR})$  and induces  an action on $N$ which preserves  the above Sasakian structure. (In fact, $N$ is isomorphic with $PSL(2,{\bR})$  and via this isomorphism  the above action is precisely the usual left action of a Lie group on itself.) If now  $\Gamma \subset PSL(2, {\bR})$ is a Fuchsian group with a compact fundamental domain then  $N/\Gamma$ is a compact Seifert manifold with a natural Sasakian structure.}
\end{remark}

\subsection{Geometric Seifert structures}

The main result of this subsection shows  that any   compact, oriented, Seifert
   3-manifold admits a $(K, \lambda)$ structure. This  will follow easily  from
    the description of the {\em geometric Seifert structures} in [JN] or [S].

\begin{theorem}{\rm Any compact, oriented, Seifert 3-manifold admits a $(K, \lambda)$ structure.}
\label{th: seifert}
\end{theorem}

\noindent{\bf Proof}\hspace{.3cm} We will begin by reviewing  the basic facts
 about the geometric  Seifert structures in a form  suitable  to the application
  we  have in mind. For details we refer to [JN] or [S] and the references therein.

A geometric structure on a manifold $M$ is a complete locally homogeneous Riemann metric of finite volume.  The universal cover of a manifold $M$ equipped with a geometric structure is a  homogeneous space  which we will call the {\em  model} of the structure. It is known that if a 3-manifold admits a geometric structure then its  model 
 belongs to a list of 8   homogeneous spaces (see[S]).

Any Seifert manifold admits a geometric structure  corresponding to one of the following   6 models:
\[
S^2\times {\bE}^1,\;\;\;{\bE}^3,\;\;\;{\bH}^2 \times{\bE}^1,\;\;\;S^3,\;\;\;N,\;\;\;\tilde{PSL}
\]
where ${\bE}^k$ denotes the $k$-dimensional Euclidean space, ${\bH}^2$ denotes the hyperbolic plane, $N$ denotes the Heisenberg group equipped with a left invariant metric and $\tilde{PSL}$ denotes the universal cover of $PSL(2,{\bR})\cong {\rm Isom}\,({\bH}^2)$. 

 According to [RV], the Seifert manifolds which admit $N$ as a model are  the  nontrivial $S^1$ bundles over a torus and as we have seen in the previous subsection such manifolds admit $(K,\lambda)$ structures.

The Seifert manifolds  which admit  geometric structures  by  ${\bE}^3$ are flat space form and are completely described in [Wo]. One can verify directly that these admit  natural $(K, \lambda)$ structures.

$S^3$ has a natural $(K,1)$ structure as the total space of the Hopf fibration $S^3 \ra S^2$. Any Seifert manifold modeled by $S^3$ is obtained as the quotient by a finite group of  {\em fiber preserving isometries}. Thus they all inherit a $(K,1)$ structure.  

If ${\bX}$ is a model other than $S^3$ or ${\bE}^3$ then the group of isomorphisms which fix a given point $x\in {\bX}$ fixes a tangent direction at that point. So ${\bX}$ has an ${\rm Aut}\,({\bX})$-invariant tangent line field. This line field  fibers  ${\bX}$ over ${\bS}^2$, ${\bE}^2$ or ${\bH}^2$.

 For ${\bX}= {\bS}^2 \times {\bE}^1, {\bH}^2 \times {\bE}^1$ this is the obvious fibration.   $\tilde{PSL}$ can be alternatively identified  with the universal cover of the unit tangent bundle  $T_1{\bH}^2$ of ${\bH}^2$. It thus has a natural line fibration which coincides with the  fibration abstractly described above. Note that the $(K, 1)$ structure on $T_1{\bH}_2$ constructed at the end of $\S 1.3$ lifts to the universal cover $\tilde{PSL}$.

If ${\bX}$ is one of  of these remaining three models  denote by ${\rm Aut}_f(X)\subset {\rm Aut}\,(X)$ the  subgroup  preserving the above line fibration  (as an oriented fibration). Note that each of them admits a $(K, \lambda)$ structure and ${\rm Aut}_f$  is in fact a group of isomorphisms of this structure.

 The trivial Seifert manifold  $S^2 \times S^1$   presents a few ``pathologies'' as far as geometric Seifert structures are concerned  (see [JN]) but we do not need to worry since it obviously admits a $(K,0)$ structure.

The other Seifert manifolds   which admit geometric structures modeled by $S^2 \times {\bE}^1$, ${\bH}^2 \times {\bE}^1$ or $\tilde{PSL}$ can be obtained as quotients $\Pi \setminus  {\bX}$ where $\Pi$ is some subgroup  of ${\rm Aut}_f$. Thus $\Pi$ invariates the  universal $(K,\lambda)$ structures on  these models and therefore the quotients  will admit such structures 
as well. 

The list of Seifert manifolds is complete. Note in particular that the Seifert manifolds  geometrized by $N$,  ${\bH}^2 \times {\bE}^1$ and $\tilde{PSL}$ are $K(\pi,1)$'s and hence they cannot admit  metrics of  nonnegative scalar curvature. \hspace{.3cm} $\Box$

\begin{remark}{\rm The above analysis can be refined to offer an answer to the 
question raised in [We]: which Seifert manifolds admit Sasakian structures. The 
answer is simple. A Seifert manifold admits  a Sasakian structure if and only if
 its (rational) Euler class is negative.  According to [NR], these are precisely
  the Seifert manifolds which can occur as links  of a quasi-homogeneous singularity.
   This extends (in the $3D$  case)  the previous result of [Sas] concerning Sasakian 
   structures on Brieskorn manifolds.

This fact was observed by many other authors (see [ENS], [Ne], [V]). In fact this 
geometry  of the link is in most cases a complete invariant of the analytic
type of the singularity (see [Ne], [Sch]).}
\end{remark}

\section{Dirac operators on 3-manifolds}
\setcounter{equation}{0}

In this section we discuss  the relationships between  $spin^c$
structures  and m.a.c. structures on a 3-manifold.  On any m.a.c. 3-manifold, besides  ${\rm spin}^c$ Dirac operators  there exists another natural Dirac operator which imitates the Hodge-Dolbeault operator on a complex manifold.  We will analyze the relationships between them.

\subsection{3-dimensional spinorial algebra}

We include here a brief survey of  the basic facts about the representations of $Spin(3)\cong SU(2)$.   Denote by $Cl_3$ the Clifford algebra generated by   $V={\bR}^3$ and by ${\bH}$ the skew-field of quaternions. Consider an orthogonal basis    $\{e_0,e_1,e_2\}$ of  $V$. 

It is convenient to identify  ${\bH}$ with ${\bC}^2$ via the correspondence
\be
{\bH} \ni q = u+ {\bf j}v  \mapsto \left[
\begin{array}{r}
u \\
v
\end{array} \right] \in {\bC}^2
\label{eq:  qua}
\ee
where
\[
q=a+b{\bf i}+ c{\bf j} + d{\bf k},\;\;u=a+b{\bf i},\;\;v =c-d{\bf i}.
\]
For each quaternion $q$ denote by $L_q$  (resp. $R_q$) the left (resp. the right) multiplication by $q$.
$R_{\bf i}$ defines a complex structure on ${\bH}$ and the correspondence (\ref{eq: qua}) defines an isomorphism of complex vector spaces.

The Clifford algebra $Cl_3$ can be represented on ${\bC}^2\cong {\bH}$  using the correspondences
\be
e_0 \mapsto L_{\bf i}  \longleftrightarrow {\bf c}(e_0)= \left[
\begin{array}{rr}
{\bf i}  & 0 \\
0 & -{\bf i} \end{array}
\right].
\label{eq: qua0}
\ee
\be
e_1 \mapsto L_{\bf j} \longleftrightarrow {\bf c}(e_1)= \left[
\begin{array}{rr}
0  & 1 \\
-1 & 0 \end{array}
\right].
\label{eq: qua1}
\ee
\be
e_2 \mapsto L_{\bf k} \longleftrightarrow {\bf c}(e_2)=\left[
\begin{array}{rr}
0&{\bf i}  \\
 {\bf i}&0 \end{array}
\right].
\label{eq: qua2}
\ee
The restriction of  the above representation to   $Spin (3)\subset Cl_3$ defines the   complex spinor representation of  $Spin(3)$
\[
{\bf c}: Spin(3) \ra {\rm Aut}\, ({\bS}_3).
\]
The Clifford multiplication  map 
\be
{\bf c}:\Lambda^*V \stackrel{\gq}{\ra} Cl_3\ra  {\rm End}\,({\bS}_3)
\label{eq: correspond}
\ee
 identifies $\Lambda^1V$ with the space of traceless, skew-hermitian endomorphisms of ${\bS}_3$. It extends by complex linearity to a map from $\Lambda^1V \otimes {\bC}$ to the space  of traceless endomorphisms of ${\bS}_3$. In particular, the purely imaginary
1-forms are mapped to selfadjoint endomorphisms.

For each $\phi \in {\bS}_3$ denote by  $\tau(\phi)$ the endomorphism of  ${\bS}_3$ defined by
\[
\tau(\phi) \psi =\lan \psi , \phi \ran \phi - \frac {1}{2}|\phi|^2 \psi.
\]
The map $\tau$ plays a central role in the 3-dimensional Seiberg-Witten equations. 

If we identify ${\bS}_3$  with  ${\bC}^2$ as above  and if $\phi$ has the form
\[
\phi =\left[
\begin{array}{r}
\alpha \\
\beta
\end{array}
\right]
\]
the $\tau(\phi)$ has the description
\be
\tau(\phi)=\frac{1}{2}\sum_{i=0}^2\lan\phi, {\bf c}(e_i)\phi\ran {\bf c}(e_i)=\left[
\begin{array}{lr}
\frac{1}{2}(|\alpha|^2 -|\beta|^2) & \alpha\overline{\beta} \\
\overline{\alpha} \beta &  \frac{1}{2}(|\beta|^2-|\alpha|^2) 
\end{array}
\right].
\label{eq: quadr}
\ee
It is not difficult to check  that the nonlinear map
\[
\tau: {\bS}_3 : \ra {\rm End}\,({\bS}_3)
\]
is $Spin(3)$ equivariant.

We want to  describe some of the invariant-theoretic features of the structure:
\[
\mbox{ (oriented Euclidean  3 dimensional space + distinguished unit vector)}. 
\]
This is the algebraic counterpart of  a m.a.c. structure on an oriented 3-manifold.

Assume $V$ is an oriented  Euclidean space which has a distinguished unit vector, say $e_0$. The  group  of isomorphisms of  this structure is $U(1)\cong S^1\cong SO(2)$. The group $Spin(3)$  acts naturally on $V$. The Lie algebra of the subgroup $H$ of $Spin(3)$ which fixes $e_0$  is generated by $e_1e_2 ={\goq}(\ast e_0)$ and can be  identified with ${\uu}(1)$    via the correspondence
\[
e_1e_2 \mapsto {\bf i} \in {\uu}(1).
\]
This tautologically identifies $H$ with $S^1$.  The representation of $H$ on ${\bS}_3$ is no longer irreducible and  consequently ${\bS}_3$ splits as a direct  sum of irreducible $H$-modules.  Alternatively, this splitting can be described  as the unitary spectral decomposition of ${\bS}_3$ defined by the action of $e_1e_2$ on ${\bS}_3$. According to (\ref{eq:
qua1}) and (\ref{eq: qua2}) we have 
\[
{\bS}_3 \cong {\bS}_3({\bf i})\oplus {\bS}_3(-{\bf i}).
\]
The action of $H$ on ${\bS}_3({\bf i})$ is the tautological $S^1$ representation, while  the action on ${\bS}_3(-{\bf i})$ is the conjugate of the tautological representation.

\subsection{3-dimensional spin geometry}

Consider a compact, oriented,  m.a.c. 3-manifold $(N,\eta,g)$. Since $w_2(N)=0$
the manifold $N$ is spin.   To understand the relationship between the spin structures on $N$ and  the m.a.c structure we need to consider  gluing data of $TN$ compatible with the m.a.c. structure.

Consider a  good, open cover $\{U_\alpha\}$ of $N$ and a gluing cocycle 
\[
g_{\alpha \beta}: U_{\alpha \beta} \ra SO(2)\cong U(1)
\] 
defining $TN$. The cocycle  is valued in $SO(2)$ since $TN$ has a distinguished section $\zeta$, the dual of $\eta$.  Note that $g_{\alpha \beta}$ defines a complex structure in the real 2-plane bundle $\lan \zeta\ran^\perp$. It is not difficult to see that  $\lan \zeta \ran^\perp \cong {\can}^{-1}$ {\em as complex line bundles.}

A spin  structure on $TN$ is a lift of  this cocycle to an $H$-valued cocycle,
 where $H$ is the subgroup of $Spin(3)$ defined at the end of  {\S}2.1. Using 
 the tautological identification $H\cong S^1$ we can identify the cover $H \ra SO(2)$ 
  with 
\[
S^1 \stackrel{z^2}{\ra}S^1.
\]
Hence a spin structure is defined by a cocycle 
\[
\tilde{g}_{\alpha \beta}:U(1)\cong SO(2) \ra U(1)
\]
such that $\tilde{g}_{\alpha \beta}^2 =g_{\alpha \beta}$. In other words, $\tilde{g}_{\alpha \beta}$ defines a  square root of ${\can}^{-1}$. 
Moreover, two such lifts define isomorphic square roots if and only if they are cohomologous so  there  exists a bijective correspondence between the square roots of ${\can}^{-1}$  (or, equivalently ${\can}$) and the  spin structures on $N$.

Now, fix a spin structure  on $N$ defined by a lift
\be
\tilde{g}_{\alpha \beta}: U_{\alpha \beta} \ra H.
\label{eq: princ}
\ee
 The  complex  spinor  bundle ${\bS}$  of this spin structure is  associated
 to the principal $H$-bundle defined by (\ref{eq: princ}) via the
 representation
 \[
  H \hookrightarrow Spin(3) \ra {\rm Aut}\,({\bS}).
 \]
 As we have already seen this splits as $\tau_1 \oplus  \tau_{-1}$ where $\tau_1$ denotes the tautological  representation of $S^1\cong H$ and 
 $\tau_{-1}$ is its conjugate.
 
 The component $\tau_1$ defines the  square root of of $\lan \zeta \ran^\perp$
 i.e. the  line bundle ${\can}^{-1/2}$ characterizing  the  chosen spin
 structure.   We have thus shown  that a choice of a m.a.c. structure on $N$
 produces a splitting
 \be
 {\bS} \cong {\can}^{-1/2}\oplus {\can}^{1/2}.
 \label{eq: split}
 \ee
 
Once we have fixed a spin structure it is very easy to classify the ${\rm
spin}^c$  structures:  they are bijectively parameterized by the complex line
bundles $L\ra N$.  The complex  spinor bundle associated to the ${\rm spin}^c$
structure defined  by the line bundle $L$ is
\[
{\bS}_L \cong {\bS}\otimes L \cong {\can}^{-1/2} \otimes L  \oplus
{\can}^{1/2}\otimes L
\]
An important  important special case is when $L = {\can}^{-1/2}$. In this case
\[
{\bS}_{{\can}^{-1/2}} \cong  
{\can}^{-1}\oplus {\uc}\stackrel{def}{=}{\bS}_\eta
\]
where we denoted by ${\uc}$ the trivial complex line bundle over $N$.

\begin{remark}{\rm  Our sign conventions differ from those of  [MST].  There they chose ${\bS}({\bf i}) \cong {\can}^{1/2}$. The overall effect is a permutation of rows  and columns in the  block description of the geometric Dirac operator of $\S 2.4$.}
\end{remark}

\subsection{Pseudo Dolbeault operators}

The complex bundle ${\bS}_\eta$ introduced in $\S 2.2$ was {\em a priori} defined in terms of
 a fixed spin structure on $N$ but {\em a posteriori}, the spin structure becomes
irrelevant. This is similar to complex  manifolds  where $\Lambda^{0,*}T^*$  is
a complex spinor bundle of  the  ${\rm spin}^c$ structure
canonically associated to the complex  manifold.  In that case the Dolbeault
operator is a Dirac operator compatible with the Clifford structure. Moreover
it is a geometric Dirac operator if the manifold is K\"{a}hler.

In this subsection want to pursue this analogy a little further.  In the
process we will construct an operator which behaves very much like the
Dolbeault operator in the complex case.

Denote by ${\cal P}$ the  real 2-plane bundle $\lan\eta \ran^\perp$. We orient ${\cal P}$ using the complex structure $-{\gc}(\ast \eta)$ which identifies it with ${\can}$.   Consider now a complex hermitian vector bundle $E\ra N$. Any connection  $\nabla$ on $E$ defines an operator
\[
\nabla  :C^{\infty}(E)\ra C^{\infty}(T^*N\otimes  E).
\]
Now observe that
\[
T^*N\otimes  E \cong \left(\lan \eta \ran  \oplus {\cal P}\right) \otimes E \cong E
\, \oplus \, ({\can}\otimes E )\, \oplus \,({\can}^{-1}\otimes E).
\]
Hence for any section $\psi \in C^{\infty}(E)$ the covariant derivative $\psi$ orthogonally splits  into three components:
\[
\nabla_\zeta \psi \in C^{\infty}(E)
\]
\[
\naf \psi \in C^{\infty}({\can}\otimes E)
\]
and
\[
\onaf \psi \in C^{\infty}({\can}^{-1} \otimes E).
\]
It terms of a local adapted frame $\{\zeta_0=\zeta, \zeta_1, \zeta_2\}$ we have
\[
\naf \psi ={\ve}\otimes (\nabla_1 -{\bf i} \nabla_2)\psi
\]
\[
\onaf\psi=\bar{\ve}\otimes (\nabla_1 +{\bf i}\nabla_2)
\]
where
\[
\nabla_j =\nabla_{\zeta_j},\;\;{\ve}=\frac{1}{\sqrt{2}}(\eta^1+{\bf i}\eta^2),\;\;\bar{\ve}=\frac{1}{\sqrt{2}}(\eta^1-{\bf i}\eta^2).
\]
For example when $E={\uc}$ then $\onaf \psi \in C^{\infty}({\can}^{-1})$. In this case if $\nabla$ is the trivial connection $d$ we will write $\partial_\flat$ (resp. $\bar{\partial}_\flat$) instead of $\naf$ (resp. $\onaf$). 
 Notice that if  $E= {\can}^{-1}$ then $\naf \psi \in C^{\infty}({\uc})$.
 
 Coupling a  metric connection  $\nabla^E$ on $E$ with the Levi-Civita
 connection $\nabla^\perp$ on ${\can}^{\pm}$ we obtain connections
 $\nabla^{E, \pm}$ on $E\otimes {\can}^{\pm}$ and in particular operators
 \[
 \naf^{E,-}:C^\infty(E\otimes {\can}^-)\ra
 C^\infty(E),\;\;\onaf^{E,+}:C^\infty(E\otimes {\can})\ra C^\infty(E).
 \]
 Using the above explicit form of the operators $\naf$ and $\onaf$  and the
 structural equations of the background m.a.c structure we obtain the following
 result.
 
 \begin{lemma}{\rm For any hermitian vector bundle $E$ and any hermitian
 connection $\nabla^E$ on $E$ we have
 \[
 \naf^*=\onaf^{E,+},\;\;\onaf^*=\naf^{E,-}
 \]
 where the upper $*$ denotes the formal adjoint.}
 \label{lemma: formadj}
 \end{lemma}

The Levi-Civita connection on $(N,g)$ induces via orthogonal projection a connection $\nabla^\perp$ on ${\can}^{-1}$ compatible with the hermitian structure. The {\em pseudo Dolbeault} operator on $(N,\eta, g)$ is the first order partial differential operator ${\dol}:C^{\infty}({\bS}_\eta) \ra C^{\infty}({\bS}_\eta)$ which in terms of the splitting ${\bS}_\eta\cong {\can}^{-1}\otimes {\uc}$ has the block decomposition
\[
{\dol}_N=\left[
\begin{array}{rr}
{\bf  i}\nabla^\perp_\zeta  & \bar{\partial}_\flat \\

 &  \\

\naf^\perp  & -{\bf i}\partial_\zeta 
\end{array}
\right]
\]
 where $\naf^\perp$ is obtained as above starting from the connection
 $\nabla^\perp$ on ${\can}^{-1}$. More generally, consider the complex spinor bundle ${\bS}_L$ associated to the ${\rm spin}^c$ structure defined by the complex line bundle $L$
\[
{\bS}_L\cong {\can}^{-1/2}\otimes L  \oplus {\can}^{1/2}\otimes L.
\]
Using the connection $\nabla$ on ${\can}$ and  a connection $A$  on $L$  we obtain connections $\nabla^{\pm}$  on  ${\can}^{\pm 1/2} \otimes L$. We  can produce a twisted pseudo-Dolbeault operator ${\dol}_L={\dol}_{L, A}$ on ${\bS}_L$ described  by the block decomposition
\[
{\dol}_L={\dol}_{L,A}=\left[
\begin{array}{cr}
{\bf i}\nabla_\zeta^{-}  & {\onaf}^{+} \\

 &  \\

\naf^- & -{\bf i}\nabla_\zeta^{-}
\end{array}
\right].
\]
Lemma \ref{lemma: formadj} shows that ${\dol}_L$ is formally selfadjoint. Note that
\[
{\dol}_N ={\dol}_{{\can}^{-1/2}}.
\]
When $L$ is trivial we set ${\dol}_L  ={\dol}_0$.

Denote by ${\bf F}^{\pm}$ the curvature of ${\nab}^\pm$.  We want to analyze the ``commutator''
\[
[{\obnaf}^+, {\nab}^+_\zeta] \stackrel{def}{=}{\obnaf}^{+}\circ {\nab}_\zeta^+  - {\nab}_\zeta^- \circ {\obnaf}^+.
\]
assuming $N$ is a Killing m.a.c manifold.  Choose a local adapted frame $\{\zeta,\zeta_1, \zeta_2\}$ with dual coframe $\{\eta, \eta_1, \eta_2\}$. For any $\psi \in C^{\infty}({\can}^{1/2}\otimes L)$ we have
\[
[{\obnaf}^+, {\nab}^+_\zeta]\psi=\bar{\ve}\otimes({\nab}^+_1+{\bf i}{\nab}^+_2){\nab}^+_\zeta \psi -{\nab}_\zeta^-\{\bar{\ve}\otimes ({\nab}^+_1+{\bf i}{\nab}^+_2)\}\psi
\]
\[
=\bar{\ve}\otimes \left({\nab}^+_1+{\bf i}{\nab}_2^+\right){\nab}^+_\zeta \psi-({\na}^\perp\bar{\ve})\otimes \left({\nab}^+_1 +{\bf i}{\nab}^+_2\right)\psi
\]
\[
-\bar{\ve}\otimes \left({\nab}_\zeta^+({\nab}^+_1 +{\bf i}{\nab}^+_2)\right)\psi
\]
\[
=\bar{\ve}\otimes\left\{({\nab}^+_1 +{\bf i}{\nab}^+_2){\nab}^+_\zeta -{\nab}^+_\zeta({\nab}^+_1+{\bf i}{\nab}^+_2)\right\}\psi -({\na}^\perp_\zeta
\bar{\ve})\otimes({\nab}^+_1+{\nab}^+_2)\psi
\]
(use $\nabla^\perp_\zeta \bar{\ve}=-{\bf i}{\vfi}(x)\bar{\ve}$)
\[
=\bar{\ve}\otimes( {\bf F}^+(\zeta_1, \zeta) +{\bf i}{\bf F}^+(\zeta_2, \zeta) +{\nab}^+_{[\zeta_1,\zeta]}+{\bf i}{\nab}^+_{[\zeta_2,\zeta]})\psi
+ {\bf i}{\vfi}(x) \bar{\ve}\otimes({\nab}^+_1 +{\bf i}{\nab}^+_2)\psi
\]
(use (\ref{eq: commu1}) and (\ref{eq: commu2}) )
\[
=\bar{\ve}\otimes({\bf F}^+(\zeta_1, \zeta) +{\bf i}{\bf F}^+(\zeta_2, \zeta)  +b(x){\nab}^+_2 -{\bf i}b(x){\nab}^+_1)\psi 
+ {\bf i}{\vfi}(x) \bar{\ve}\otimes({\nab}^+_1 +{\bf i}{\nab}^+_2)\psi
\]
\[
=\bar{\ve}\otimes\left\{ {\bf F}^+(\zeta_1, \zeta) +{\bf i}{\bf F}^+(\zeta_2, \zeta)  -{\bf i} b(x)({\nab}^+_1 +{\bf i}{\nab}^+_2)\right\}\psi
+{\bf i}{\vfi}(x) \bar{\ve}\otimes({\nab}^+_1 +{\bf i}{\nab}^+_2)\psi
\]
\[
=\bar{\ve} \otimes({\bf F}^+(\zeta_1, \zeta) +{\bf i}{\bf F}^+(\zeta_2, \zeta)  )
\psi -{\bf i} \lambda  \bar{\ve}\otimes ({\nab}^+_1+{\bf i} {\nab}^+_2 )\psi.
\]
We have thus proved
\be
[{\obnaf}^+, {\nab}^+_\zeta]  =-{\bf i}\lambda(x)\, {\obnaf}^+  + \bar{\ve}\otimes\left\{{\bf F}^+(\zeta_1, \zeta) +{\bf i}{\bf F}^+(\zeta_2, \zeta)  \right\}
\label{eq: swcommute} 
\ee

When dealing with the Seiberg-Witten equations it is convenient  to describe the curvature term in the above formula in terms of the curvature $F_A$ of  $L$.   We will use the formula
\[
{\bf F}^+ =F_{{\can}^{1/2}\otimes L}= F_A +\frac{1}{2}F_{\can}.
\]
Hence we need to explicitly describe the  curvature of ${\can}$ equipped with the connection $\nabla^\perp$.

 Note that ${\can}^{-1}$ can be identified with the bundle $\lan \zeta \ran^\perp$ equipped with the complex structure
\[
{\bf i} \zeta_1  =\zeta_2\;\;\;{\bf i}\zeta_2 = -\zeta_1.
\]
 We will compute the curvature of this line bundle using  the structural equations of $\nabla^\perp$
\[
\nabla^\perp \zeta_j = -{\bf i}C \otimes \zeta_j \;\; j=1,2.
\]
Assuming $N$ is a Killing  m.a.c manifold we deduce after  some simple manipulations
\[
F_{{\can}^{-1}}(\zeta_1, \zeta)=-{\bf i}(-\partial_1 {\vfi}(x) +\partial_\zeta C_1 + b(x) C_2)
\]
\[
F_{{\can}^{-1}}(\zeta_2, \zeta)=-{\bf i}(-\partial_2 {\vfi}(x) +\partial_\zeta C_2 -b(x) C_1).
\]
Set $\mu(x) = C_1 +{\bf i} C_2$.   Some elementary algebra  shows 
\[
F_{{\can}^{-1}}(\zeta_1 , \zeta) + {\bf i}F_{{\can}^{-1}}(\zeta_2, \zeta)= {\bf i}(\partial_1 +{\bf  i}\partial_2){\vfi}
- ({\bf i}\partial_\zeta +b(x))\mu.
\]
Hence
\be
\bar{\ve}\otimes (F_{\can}(\zeta_1, \zeta)+{\bf i}F_{\can}(\zeta_2, \zeta) = -{\bf  i} \bar{\partial}_\flat {\vfi} +\bar{\ve}\otimes\{ ({\bf i}\partial_\zeta +b(x))\mu \}.
\label{eq: cancurv}
\ee
We now want to clarify the ``mysterious'' term $({\bf i}\partial_\zeta +b(x))\mu $ in the above formula.  To achieve this we will use the structural equations (\ref{eq: structural}). Thus
\[
d\eta^1 = A\wedge \eta - C \wedge \eta^2  =-b(x)\eta \wedge \eta^2 -C_1 \ast \eta
\]
and
\[
d\eta^2 =-B\wedge \eta + C\wedge \eta^1 =b(x) \eta \wedge \eta^1 -C_2\ast \eta.
\]
Temporarily set $\omega =\sqrt{2}{\ve} =\eta^1  +{\bf i} \eta^2$. The above equalities  yield
\be
d\omega = {\bf i} b(x) \eta \wedge \omega  -\mu \ast \eta.
\label{eq: myster}
\ee
Differentiating the last equality we deduce
\[
d(\mu \ast \eta) ={\bf i} d(b(x)  \eta \wedge \omega)
\]
i.e.
\[
(\partial_\zeta \mu) dvol_g \,= \, {\bf i}\{ db(x) \wedge \eta \wedge \omega + b(x) d\eta \wedge \omega - b(x) \eta \wedge d\omega\}.
\]
The middle term in the right-hand-side of the above formula cancels. The third term can be  computed using (\ref{eq: myster}). Hence
\[
(\partial_\zeta \mu) dvol_g ={\bf i}\{db(x) \wedge \eta \wedge \omega +b(x) \eta \wedge  \ast \eta\}
\]
or equivalently,
\[
(\partial_\zeta -{\bf i} b(x))\mu \,dvol_g \,  = \,{\bf i} d b(x) \wedge \eta \wedge \omega.
\]
Since
\[
db(x) \wedge \eta \wedge \omega  = (-{\bf i} \partial_1b(x) +\partial_2 b(x) )d vol_g
\]
we deduce
\[
({\bf i} \partial_\zeta + b(x) )\mu ={\bf i} (\partial_1 +{\bf i} \partial_2)b(x).
\]
In a more invariant form
\be
\bar{\ve}\otimes ({\bf i} \partial_\zeta + b(x) )\mu  = {\bf i}\bar{\partial}_\flat b(x).
\label{eq: myster1}
\ee
Using  the above equality in  (\ref{eq: cancurv}) we deduce
\be
\bar{\ve}\otimes (F_{\can}(\zeta_1, \zeta)+{\bf i}F_{\can}(\zeta_2, \zeta)\} =-{\bf  i} \bar{\partial}_\flat {\vfi} + {\bf i}\bar{\partial}_\flat b(x ) = {\bf i} \bar{\partial}_\flat \lambda(x).
\label{eq: holocurv}
\ee
Set
\[
F_A^{0,1}:=\bar{\ve}\otimes F_A(\zeta, \zeta_1+{\bf i}\zeta_2)
\]
and 
\[
F_A^{1,0}:={\ve}\otimes F_A(\zeta, \zeta_1 -{\bf i}\zeta_2).
\]
We can now rephrase the commutativity relation (\ref{eq: swcommute}) as
\be
[{\obnaf}^+, {\nab}^+_\zeta]  =-{\bf i}\lambda(x)\, {\obnaf}^+   + \frac{\bf i}{2} \bar{\partial}_\flat \lambda(x) - F_A^{0,1}
\label{eq: excellent}
\ee
By passing to  formal adjoints  we deduce  using Lemma \ref{lemma: formadj}
\be
[\naf^-, \nabla_\zeta^-]=[({\obnaf}^+)^*, {\nab}^-_\zeta]:=({\obnaf}^+)^*{\nab}^-_\zeta -{\nab}^+_\zeta ({\obnaf}^+)^*= {\bf i}\lambda ({\obnaf}^+)^* -\frac{\bf i}{2}({}^\flat\bar{\partial})^* \lambda -F^{1,0}_A.
\label{eq: excellent1}
\ee

On a $(K,\lambda)$ manifold the commutativity relations (\ref{eq: excellent}) and (\ref{eq: excellent}) further simplify.  Fix a ${\rm spin}^c$ structure
determined by a line bundle $L$  and choose a connection $A$ on $L$.  It will be extremely convenient to introduce two new differential operators  $Z_A, T_A:C^{\infty}({\bS}_L)\ra C^\infty({\bS}_L)$ defined by the block  decompositions
\[
Z_A=\left[\begin{array}{cc}
{\bf i}\nabla_\zeta^{A,-}  & 0 \\
0   & -{\bf i} \nabla_\zeta^{A, +} 
\end{array}\right]
\]
\[
T_A=\left[ \begin{array}{cc}
0  & {\obnaf}^+ \\
({\obnaf}^+)^*  & 0
\end{array} 
\right]
\]
Then the commutativity relations  (\ref{eq: excellent}) and (\ref{eq: excellent1}) can be  simultaneously rephrased as an anti-commutator identity
\be
\{Z_A, T_A\} := Z_AT_A +T_A Z_A =  -\lambda T_A -{\bf i}\left[\begin{array}{cc}
0 & -F_A^{0,1} \\
F_A^{1,0}  & 0 
\end{array}
\right].
\label{eq: zt}
\ee

We want to emphasize here that the operators $Z_A$, $T_A$ and $F_A^{0,1}$
{\em depend on the background m.a.c. structure}. Both $Z_A$ and $F_A^{0,1}$
{\em will be affected by anisotropic deformations} while $T_A$  is invariant.

\subsection{Geometric Dirac operators}

In this subsection we  will analyze the geometric Dirac operators on a 3-manifolds and in particular we will relate them with the  pseudo-Dolbeault operators of $\S 2.3$.

Consider an oriented  Killing m.a.c manifold $(N,\eta, g)$ with a {\em fixed} spin structure. Denote by ${\bS}\cong {\can}^{-1/2}\oplus {\can}^{1/2}$ the bundle of complex spinors associated to this structure.

We begin by first recalling the construction of  the canonical connection on ${\can}$. Pick a local adapted frame
$\{\zeta_0=\zeta,\zeta_1, \zeta_2\}$ and denote by $\sigma_j$ the Clifford multiplication by $\zeta_j$, $j=0,1,2$.   With respect to the canonical decomposition ${\bS}\cong {\can}^{-1/2}\oplus {\can}^{1/2}$ these operators have the descriptions 
\[
\sigma_0 =\left[
\begin{array}{rl}
{\bf i} & 0 \\
0 & -{\bf i} \end{array}
\right]
\]
\[
\sigma_1= \left[
\begin{array}{ll}
0 & \bar{\ve} \\
-{\ve } & 0 
\end{array}
\right]
\]
\[
\sigma_2 =\left[
\begin{array}{rr}
0 & {\bf i}\bar{\ve} \\
{\bf i}{\ve} & 0 
\end{array}
\right]
\]
where ${\ve}$ (resp. $\bar{\ve}$) denotes  the tensor multiplication by $\bar{\ve}$ (resp ${\ve}$)
\[
\bar{\ve}:{\can}^{1/2}\ra {\can}^{-1/2}\;\;\; ({\rm resp.}\;\;\;{\ve}:{\can}^{-1/2}\ra {\can}^{1/2}).
\]
If  $(\omega_{ij})$ denotes the ${\uso}(3)$ valued 1-form associated to the Levi-Civita connection via the local frame $\{\zeta_j\}$ i.e.
\[
\nabla \zeta_j =\sum_i \omega_{ij}\zeta_i
\]
then the canonical connection on ${\bS}$ is defined by
\[
\nabla = d -\frac{1}{2}\sum_{i<j}\omega_{ij}\otimes \sigma_i \sigma_j.
\]
Using the structural equations (\ref{eq: levicivita}) we deduce that, with respect to the local frame $\{\zeta_j\}$,  the canonical connection  on ${\bS}$ has the form
\[
{\nah}=\nabla^{\bS}= d -\frac{1}{2}(A\otimes \sigma_2 + B\otimes \sigma_1 + C\otimes \sigma_0).
\]
Using the fact that $N$ is a Killing m.a.c. we deduce
\[
{\nah}_\zeta=\partial_\zeta -\frac{1}{2}\left[
\begin{array}{rl}
{\bf i}{\vfi}(x) & 0\\
0& -{\bf i}{\vfi}(x) 
\end{array}
\right]
\]
\[
{\nah}_1=\partial_{\zeta_1}-\frac{1}{2}\left[
\begin{array}{rr}
{\bf i}C_1 & \lambda \bar{\ve} \\
-\lambda {\ve} & -{\bf i} C_1
\end{array}
\right]
\]
\[
{\nah}_2 =\partial_{\zeta_2}-\frac{1}{2} \left[
\begin{array}{rr}
{\bf i}C_2 & {\bf i}\lambda \bar{\ve} \\
{\bf i}\lambda {\ve} & -{\bf i}C_2 
\end{array}
\right].
\]
The canonical, untwisted  (geometric) Dirac operator   on ${\bS}$ is defined by 
\[
{\dir}_0= {\dir}_{\bS}=\sigma_0 {\nah}_0 +\sigma_1 {\nah}_1 +\sigma_2 {\nah}_2
\]
\[
=
\left[
\begin{array}{cc}
{\bf i}(\partial_\zeta -{\bf i}{\vfi}/2)  & \bar{\ve}\otimes\{(\partial_1 + {\bf i}C_1/2 )+{\bf i}(\partial_2+{\bf i} C_2/2)\} \\

 &  \\

{\ve}\otimes \{-(\partial_1 +{\bf i}C_1/2) +{\bf i}(\partial_2 +{\bf i}C_2/2)\} & -{\bf i}(\partial_\zeta +{\bf i}{\vfi}/2)
\end{array}
\right] +\lambda {\bf 1}_{\bS}.
\]
In terms of the pseudo-Dolbeault operator we have 
\[
{\dir}_{\bS}={\dol}_{\bS} + \lambda.
\]
More generally if we twist ${\bS}$ by a line bundle $L$ equipped with a connection $\nabla^L$ we obtain a  geometric Dirac operator on ${\bS}_L ={\bS}\otimes L$ and as above one  establishes the following identity
\be 
{\dir}_L =  {\dol}_L +\lambda.
\label{eq: dirdolbo}
\ee

\section{The Seiberg-Witten equations}
\setcounter{equation}{0}

In this section we finally  take-up the  study of the  3-dimensional Seiberg-Witten equations. We will    restrict our considerations to the special case  when the 3-manifold $N$ has  a$(K,\lambda)$-structure.
    When $\lambda=0$  it was observed by many authors  that these equations can be solved quite explicitly. The situation is   more complicated when $\lambda \neq 0$ for the reasons explained in the introduction.   We subject $N$ to  an anisotropic adiabatic deformation  so that in the limit $\lambda_\delta \ra 0$ and  study the behavior of the solutions of the Seiberg-Witten equations  as the metric degenerates.   Our study is sufficiently
    accurate to  provide  many informations about  the  exact solutions for the Seiberg-Witten equations corresponding to large $\delta$ and a certain range of
${\rm spin}^c$ structures.

\subsection{Generalities}

The goal of this subsection is to describe  the 3-dimensional Seiberg-Witten 
equations  and then derive a few  elementary consequences.

Consider $(N,g)$ a compact, oriented  m.a.c 3-manifold.  Fix a spin structure on $N$ defined by the square root ${\can}^{-1/2}$. The data entering the Seiberg-Witten equations are the  following.

\noindent {\bf  (a)} A  ${\rm spin}^c$ structure  determined by the line bundle $L$.

\noindent {\bf (b)} A connection $A$ of $L\ra N$.

\noindent {\bf (c)} A    spinor $\phi$ i.e. a section of the complex spinor bundle ${\bS}_{L}$ associated to the given  ${\rm spin}^c$ structure.

The connection $A$  defines a geometric Dirac operator  ${\dir}_A$  on ${\bS}_{L}$. The Seiberg-Witten equations are
\[
\left\{
\begin{array}{rcr}
{\dir}_{\bf A}\phi  & = & 0 \\
{\bf c}(\ast F_{\bf A})& =  & \tau (\phi)
\end{array}
\right.
\]
where $\ast$ is the Hodge $\ast$-operator of the metric $g$, $\tau$ is defined in (\ref{eq: quadr})  and  
${\bf c}$ is defined in (\ref{eq: correspond}). We will omit the symbol ${\bf
c}$ when no confusion is possible. 

The Seiberg-Witten equations have a variational nature. Fix a smooth connection
$A_0$ on $L$ and define
\[
{\gf}: L^{1,2}({\bS}_L\oplus {\bf i}T^*N) \ra {\bR}
\]
by
\[
{\gf}(\psi, a)=\frac{1}{2}\int_Na\wedge(F_{A_0}+ F_{A_0+a}) +\frac{1}{2}\int_N  \lan \psi,
{\dir}_{A_0+a}\psi \ran dv_g.
\]
\begin{lemma}{\rm The differential of ${\gf}$ at a point
${\gc}=(\phi,a)$ is}
\[
d_{{\gc}}{\gf}(\dot{\phi}, \dot{a})=\int\lan \dot{a}, {\bf c}^{-1}(\tau(\phi))-\ast F_{A_0+a}\ran
dv_g +\int_N{\re}\lan \dot{\phi}, {\dir}_{A_0+a}\phi\ran dv_g.
\]
\label{lemma: grad}
\end{lemma}

\noindent{\bf Proof}   Set $A=A_0+a$. We have
\[
d_{gc}{\gf}(\dot{\phi},\dot{a})=\frac{d}{dt}\!\mid_{t=0}{\gf}(\phi+t\dot{\phi},
 a+t\dot{a})
\]
\[
=\frac{1}{2}\int_N\dot{a} \wedge (F_{A_0}+F_{A})+\frac{1}{2}\int_N a \wedge d\dot{a}
-\int_N{\re}\lan \dot{\phi}, {\dir}_A\phi\ran dv_g -\frac{1}{2}\int_N \lan \phi,
{\bf c}(\dot{a})\phi\ran dv_g
\]
\[
{\rm (Stokes)}=\frac{1}{2}\int_N\dot{a} \wedge
(F_{A_0}+F_{A})+\frac{1}{2}\int_N\dot{a}\wedge da -\frac{1}{2}\int_N \lan \phi,
{\bf c}(\dot{a})\phi\ran dv_g - \int_N{\re}\lan \dot{\phi}, {\dir}_A\phi\ran dv_g
\]
\[
\int_N\dot{a}\wedge F_A +\frac{1}{2}\int_N \lan \phi,
{\bf c}(\dot{a})\phi\ran dv_g + \int_N{\re}\lan \dot{\phi}, {\dir}_A\phi\ran
dv_g.
\]
Since both $\dot{a}$ and $F_A$ are purely imaginary  we have 
is
\[
\int_N\lan \dot{a}, \ast F_A\ran dv_g =-\int_N \dot{a}\wedge  F_A
\]
where $\ast$ is the {\em complex linear} Hodge $\ast$-operator. On the other  
 hand, a simple computation shows that
$(\dot{a}=\sum_i\dot{a}_i\eta_i)$
\[
\int_N \lan \phi,{\bf c}(\dot{a})\phi\ran dv_g  =-\int_N \sum_i
\dot{a}_i\lan\phi, {\bf c}(\eta_i)\phi\ran dv_g 
\]
\[
=-\int_N \sum_i \lan\,\dot{a_i}\eta_i,\overline{\lan \phi, {\bf
c}(\eta_i)\phi\ran}\,\eta_i \,\ran 
dv_g =2\int_N\lan \dot{a}, {\bf c}^{-1}(\tau(\phi))\ran dv_g.
\]
Putting all the above together we get  the lemma.  $\Box$

\bigskip

The gauge group ${\gG}_L={\rm Aut}\, (L)\cong {\rm Map}\, (N,S^1)$ acts on the
space of pairs $(\psi,a)$ by
\[
\gamma\cdot(\psi,a):=(\gamma\cdot \psi, a -\gamma^{-1}d\gamma).
\]
Moreover
\[
{\gf}(\gamma\cdot (\psi, a))-{\gf}(\psi,a)=
-\int_N\gamma^{-1}d\gamma \wedge F_{A_0}=2\pi {\bf
i}\int_N\gamma^{-1}d\gamma\wedge 
c_1(A_0)
\]
Thus ${\gf}$ is unchanged by the  gauge transformations homotopic to the
constants. 
We see that the critical points of ${\gf}$ are precisely the solutions of the
Seiberg-Witten equations. In particular, the above considerations show that the
moduli space of solutions is invariant under the action of ${\gG}_L$ and so it
suffices to look at the quotient of this action.

The Seiberg-Witten equations have a more explicit description once we choose an adapted orthonormal frame $\zeta, \zeta_1, \zeta_2$.
 Using the decomposition ${\bS}_L\cong {\can}^{-1/2}\otimes L \oplus {\can}^{1/2}\otimes L$ we can  represent $\phi$ as
\[
\phi=\left[
\begin{array}{r}
\alpha
\\
\beta 
\end{array}
\right]
\]
 If $N$ is a $(K,\lambda)$- manifold and we denote by ${\nab}^\pm$ the covariant derivatives induced by $A$ on ${\can}^{\pm 1/2} \otimes L$  then  the Seiberg-Witten equations can be  rephrased as 
\be
\left\{
\begin{array}{rrrcl}
{\bf i}{\nab}^-_\zeta \alpha & +  {\obnaf}^+ \beta & +\lambda  \alpha &  = & 0 \\
 & & & & \\
({\obnaf}^+)^* \alpha & - {\bf i}{\nab}^+_\zeta \beta  & +\lambda  \beta & = &0 \\
&&&&\\

& &  \frac{1}{2}(|\alpha|^2-|\beta|^2) & = &{\bf i}F_{\nab}(\zeta_1, \zeta_2) \\
& & & & \\
& & {\bf i}\alpha \bar{\beta}&= &\bar{\ve}\otimes F_{\nab}(\zeta_1
+{\bf i} \zeta_2, \zeta) =-F_A^{0,1}
\end{array}
\right.
\label{eq: sw}
\ee
In particular, note that
\[
\eta \wedge c_1(A) = \frac{\bf i}{2 \pi}F_A(\zeta_1, \zeta_2)\eta \wedge \eta^1\wedge \eta^2 =\frac{1}{4\pi}(|\alpha|^2 -|\beta|^2) dvol_g.
\]
Hence
\be
\int_N \eta  \wedge c_1(F_{\bf A}) = \frac{1}{4 \pi}(\|\alpha\|^2-\|\beta\|^2)
\label{eq: chern}
\ee
where $\|\cdot \|$ denotes the $L^2$- norm  over $N$. 

In terms of the operators $ Z_A$ and $T_A$  defined in \S 2.3  we can rewrite the first two equations as
\[
(Z_A + T_A)\phi =-\lambda \phi.
\]
The anti-commutation relation   (\ref{eq: zt})   shows that if $(A,\phi)$ is a solution of the Seiberg-Witten  equation then
\be
(Z_A+T_A)^2\phi=Z^2_A\phi+T^2_A\phi -\lambda T_A\phi+   \left[
\begin{array}{cc} 0 & \alpha \bar{\beta}
\\
\bar{\alpha}\beta  & 0 
\end{array}
\right]\cdot \left[\begin{array}{c}
\alpha \\
\beta
\end{array}
\right].
\label{eq: zt1}
\ee

\subsection{Adiabatic limits}

We now have all the data we need  to study the behavior of the Seiberg-Witten equations as the metric is anisotropically deformed until it degenerates.  

 Let $(N.\eta, g)$  and ${\can}^{1/2}$ as above.    As usual we denote by   $(N, \eta_\delta, g_\delta)$ the anisotropic deformation  defined in $\S 1.1$.    For each $\delta\geq 1$ we  we will refer to the Seiberg-Witten equations defined  in terms of the metric $g_\delta$ as  $SW_\delta$ equations.  More explicitly these are ($|\zeta|_{g_1}\equiv 1,\;\;|\eta|_{g_1}\equiv 1$)
\be
\left\{
\begin{array}{rrrcl}
\delta{\bf i} {\nab}^-_\zeta \alpha & +  {\obnaf}^+ \beta &
+\frac{\lambda}{\delta} \alpha &  = & 0 \\
 & & & & \\
({\obnaf}^+)^* \alpha & - \delta{\bf i}{\nab}^+_\zeta \beta  &
+\frac{\lambda}{\delta}\beta & = &0 \\
&&&&\\

& &  \frac{1}{2}(|\alpha|^2-|\beta|^2)  & = & 
{\bf i}F_{\nab}(\zeta_1,\zeta_2)\\
& & & & \\
& &  {\bf i}\alpha \bar{\beta}&= &\delta\bar{\ve}\otimes F_{\nab}(\zeta_1
+{\bf i} \zeta_2, \zeta) 
\end{array}
\right.
\label{eq: swd}
\ee
The operator $Z_{A}$ depends on the metric $g_\delta$  through the Levi-Civita connection
$\nabla^\perp$ on ${\can}$   and so we will write $Z_{A,\delta}$ to emphasize this dependence.   As $\delta \ra \infty$ we have
\[
\frac{1}{\delta}Z_{A,\delta} \ra Z_{A,\infty}:=\left[
\begin{array}{cc}
{\bf i}\hat{\nabla}_\zeta^{A,-} & 0 \\
0 & -{\bf i} \hat{\nabla}_\zeta^{A,+}
\end{array}
\right]
\]
where  for each connection $A$ on $L$ we denoted by $\hat{\nabla}^{A,\pm}$ the  connection on ${\can}^{\pm 1/2}\otimes L$ obtained by tensoring the limiting connection $\nabla^{\infty}$ on ${\can}^{\pm 1/2}$ (described in  (\ref{eq: anisotrop1}))  with the connection $A$. On the other hand, $T_A$ is unaffected by the adiabatic changes in the metric.

We will fix a smooth connection $\nabla^0$ on $L$.  The  Sobolev norms will be defined in terms of this connection and its tensor products with the connections induced by the Levi-Civita connection of the {\em fixed} metric $g=g_1$.  An arbitrary  connection on $L$ will have the form
\[
\nabla^0 +A_\delta, \;\;\;A_\delta \in {\uu}(1) \otimes \Omega^1(N).
\]
For  $1\leq p \leq \infty$ we denoted by $\|\cdot \|_p$ the $L^p$-norm with respect to the metric $g=g_1$,  by $\| \cdot \|$ the $L^2$-norm with respect to the same metric while
$\|\cdot\|_\delta$ will denote the $L^2$-norm with respect to the metric
$g_\delta$).

We denote by ${\gA}_L$ the affine space of {\em smooth} connections on $L$ and by ${\cal S}_\delta$ the  collection of  gauge equivalence classes of solutions of  $SW_\delta$.

\begin{theorem}{\rm  Let  $(N, \eta, g)$ as above and fix a ${\rm spin}^c$ 
structure defined by a complex line bundle  $L$ (so that the  associated complex
 spinor  bundle has determinant $L^2$).   Assume that  for each  sufficiently
  large $\delta$ 
\[ 
{\cal  S}_\delta \neq \emptyset.
\]
Then  any sequence $\{([A_\delta], [\phi_\delta])\in {\cal S}_\delta\; ; \;  \delta  \gg 1\}$ 
admits a subsequence which converges  in  the $L^{1,2}$ topology to a pair 
\[
([A], [\phi])\in ({\gA}_L \times C^{\infty}(L))/{\rm Aut}\,(L)  
\]
satisfying the following conditions.}
\be
F_A(\zeta, \cdot)\equiv 0
\label{eq: 0}
\ee
\be
{\obnaf}^+\beta = \hat{\nabla}^{A,-}_\zeta \alpha =0
\label{eq: 1}
\ee

\be
({\obnaf}^+)^*\alpha = \hat{\nabla}^{A,+}_\zeta \beta =0
\label{eq: 2}
\ee

\be
\|\alpha\|\, \cdot \, \|\beta\| =0.
\label{eq: 3}
\ee
\label{th: vanish}
\end{theorem}

\noindent The equations (\ref{eq: 1}) and (\ref{eq: 2}) can be equivalently rephrased as
\be
Z_{A,\infty}\phi= T_A\phi =0.
\label{eq: decouple}
\ee
We can say  that the Seiberg-Witten equations decouple in the adiabatic limit.

\bigskip

\noindent {\bf Proof} \hspace{.3cm}    The proof of the theorem  relies essentially on the following  {\em uniform} estimates.
   
\begin{lemma}{\rm  There exist $R_1,\; R_2 > 0$ such that
\be
\sup |\phi_\delta(x)| \leq R_1\;\;\forall \delta \geq 1
\label{eq: sup}
\ee
and}
\be
\|F(A_\delta)\|_\infty \leq R_2 \;\;\;\forall \delta \geq 1.
\label{eq: L2}
\ee
\be
\| \,|\alpha_\delta|\, \cdot \, |\beta_\delta|\,\|_2 =O(1/\delta)\;\;\;{\rm as}\;\;\delta \ra \infty.
\label{eq: anul}
\ee
\label{lemma: uniform}
\end{lemma}

\bigskip

\noindent{\bf Proof of the lemma} As in Lemma 2 of [KM]  we deduce that
\[
\sup |\phi_\delta (x)| \leq \sup |s_\delta(x)|
\]
where $s_\delta$ denotes  the scalar curvature of $g_\delta$. The  estimate
(\ref{eq: sup}) is now a consequence of (\ref{eq: limit}) in $\S 1.2$.

To prove (\ref{eq: L2}) note first that $F(A_\delta)$ splits into two
orthogonal parts. A horizontal part
\[
F^h(A_\delta) = F_{12}(A_\delta)\eta^1 \wedge \eta^2
\]
and a vertical part
\[
F^v(A_\delta) = F_{01}(A_\delta)\eta \wedge \eta^1 +
F_{20}(A_\delta)\eta^2\wedge \eta.
\]
The third equation in (\ref{eq: swd})  coupled with  (\ref{eq: sup}) yields 
\[
\|F^h(A_\delta)\|_\infty = O(1) \;\;{\rm as}\;\;\delta \ra \infty.
\]
The fourth equality in (\ref{eq: swd}) implies
\be
\|F^v(A_\delta)\|_\infty \leq \mbox{const.} \delta^{-1} \|\, |\alpha_\delta|\cdot
|\beta_\delta|\, \|_\infty\, =O(\delta^{-1})
\label{eq: cube}
\ee
where
\[
\phi_\delta =\left[
\begin{array}{r}
\alpha_\delta  \\
\beta_\delta
\end{array}
\right].
\]
To prove (\ref{eq: anul}) we will use  (\ref{eq: zt1}) and we get  (writing $Z_{A_\delta}$ instead of the more accurate $Z_{A_\delta, \delta}$)
\[
\lambda^2_\delta \phi_\delta = (Z_{A_\delta}+T_{A_\delta})^2\phi_\delta = Z^2_{A_\delta}\phi_\delta  +T^2_{A_\delta}\phi_\delta  -\lambda_\delta T_{A_\delta} \phi_\delta   +   \left[
\begin{array}{cc}0 & \alpha_\delta \bar{\beta_\delta}
\\
\bar{\alpha_\delta}\beta_\delta  & 0 
\end{array}
\right]\left[\begin{array}{c}
\alpha_\delta\\
\beta_\delta
\end{array}
\right].
\]
Taking the inner product with $\phi_\delta$ in the above equality and then integrating by parts with respect to the metric $g_\delta$ we deduce
\[
\|Z_{A_\delta}\phi_\delta\|_\delta^2 +\|T_{A_\delta}\phi_\delta\|^2_\delta + 2\| \, |\alpha_\delta|\,\cdot \, |\beta_\delta|\, \|^2_\delta    = \lambda_\delta ^2\|\phi_\delta\|^2_\delta  + \lambda_\delta \lan  T_{A_\delta}\phi_\delta, \phi_\delta \ran_\delta.
\]  
The Cauchy inequality yields
\be
\|Z_{A_\delta}\phi_\delta\|_\delta^2 +\|T_{A_\delta}\phi_\delta\|^2_\delta +2 \| \, |\alpha_\delta|\,\cdot \, |\beta_\delta|\, \|^2_\delta    \leq  |\lambda_\delta|\cdot \|\phi_\delta\|_\delta\left(|\lambda_\delta|\cdot \|\phi_\delta\|_\delta   + \|T_{A_\delta}\phi_\delta\|_\delta \right)
\label{eq: 38}
\ee
Since $|\lambda_\delta|,\;\|\phi_\delta\|^2_\delta=O(\delta^{-1})$ we  first deduce from the above inequality that $\|T_{A_\delta}\phi_\delta\|_\delta =O(\delta^{-3/2})$ and then
\[
\|\,|\alpha_\delta|\, \cdot \, |\beta_\delta|\, \|_\delta =O(\delta^{-3/2}).
\]
Up to a rescaling this is precisely the inequality (\ref{eq: anul}). \hspace{.3cm}  $\Box$

\bigskip
 
  From the estimate (\ref{eq: L2}) we deduce  that  $A_\delta$ (modulo Coulomb gauges)  is bounded  in the $L^{1,p}$-norm for any $ 1< p < \infty$. Thus  a subsequence converges strongly in $L^p$ and weakly in $L^{1,p}$ to some H\"{o}lder continuous connection $A$ on $L$.
  The fourth  equation in (\ref{eq: swd}) and the estimate (\ref{eq: anul}) we
  have just proved show that 
  \[
  \|F_{A_\delta}^{0,1}\|= O(\delta^{-2}).
  \]
  The condition (\ref{eq: decouple}) is a consequence of the following
  auxiliary result.
  
  \begin{lemma} {\bf (Adiabatic decoupling lemma)} {\rm Consider a sequence of
  smooth connections  $A_\delta \in {\gA}_L$ satisfying the following
  conditions.
  
  \noindent (i) $A_\delta$  converges  in the weak $L^{1,p}$ topology $(p > 3=\dim
  N$) to a H\"{o}lder continuous connection $A$.
  
  \noindent  (ii) $\|F^{0,1}_{A_\delta}\|=o(\delta^{-1})$ as $\delta \ra \infty$.
  
  Then any sequence $\phi_\delta \in C^{\infty}({\bS}_L)$ such that
  \[
  \|\phi_\delta\|= O(1)\;\;{\rm as}\; \delta \ra \infty
  \]
  and
  \[
  {\dir}_{A_\delta}\phi_\delta =0
  \]
  contains a subsequence which converges strongly  in $L^{1,2}({\bS}_L)$ to  a spinor $\phi\in
  L^{1,2}({\bS}_L)$ satisfying}
  \be
  Z_{A,\infty}\phi=T_A\phi=0.
  \label{eq: adiabaticdec}
  \ee
  \end{lemma}
  
  \noindent{\bf Proof} \hspace{.3cm}   Set $Z_{A_\delta}=Z_{A_\delta,\delta}$
  and
  \[
  \Xi_\delta =\frac{1}{\delta} Z_{A_\delta}
  \]
  The coefficients of $\Xi_\delta$ converge uniformly to the  coefficients of
  $Z_{A,\infty}$. More precisely, the condition (i) implies $\Xi_\delta  =Z_{A,\infty}+ R$ where
  the zeroth order term $R$ satisfies
  \be
  \|R\|_\infty = o(1)
  \label{eq: error}
  \ee
  Note that although $\Xi_\delta$ is defined  using the metric $g_\delta$ it is
  formally self-adjoint with respect to the metric $g_1$ as well.
  
  The equation ${\dir}_{A_\delta}\phi_\delta=0$ can be rewritten as
  \[
  (Z_{A_\delta}+T_{A_\delta}) \phi_\delta =-\lambda_\delta \phi_\delta.
  \]
  Using the equation (\ref{eq: zt}) (with the  background m.a.c. structure $(g_\delta,
  \zeta_\delta=\delta\zeta)$) we get
  \be
  \lambda^2_\delta \phi_\delta = (Z_{A_\delta}^2 +T_{A_\delta})^2\phi_\delta
  =-\lambda_\delta T_{A_\delta} +{\cal F}_\delta
  \label{eq: ztf}
  \ee
  where
  \[
  {\cal F}_\delta=-{\bf i}\delta\left[\begin{array}{cc}
0 & -\bar{\ve}F_{A_\delta}(\zeta, \zeta_1+{\bf i}\zeta_2) \\
{\ve}F_{A_\delta}(\zeta, \zeta_1 -{\bf i} \zeta_2)  & 0 
\end{array}
\right]= -{\bf i}\delta\left[\begin{array}{cc}
0 & -F_{A_\delta}^{0,1} \\
F_{A_\delta}^{1,0}  & 0 
\end{array}
\right].
\]
Take the inner product with $\phi_\delta$  of both sides in (\ref{eq: ztf}) and
then integrate by parts using {\em the  fixed metric $g_1$}. We get
\[
\delta^2\|\Xi_\delta\phi_\delta^2\|^2 +\|T_{A_\delta}\phi_\delta\|^2 \leq
\|\phi_\delta\|\left\{ |\lambda_\delta|\cdot \|T_{A_\delta}\|+\|{\cal
F}_\delta\|+ \lambda^2_\delta \|\phi_\delta\|\right\}.
\]
Since $\|\phi_\delta\|=O(1)$ and $\|{\cal F}_\delta\|=O(\delta
\|F_{A_\delta}^{0,1}\|)= o(1)$ we deduce
\[
\delta^2\|\Xi_\delta\phi_\delta^2\|^2 +\|T_{A_\delta}\phi_\delta\|^2 \leq
C\left(\delta^{-1}\|T_{A_\delta}\phi_\delta\| +\delta^{-2} + o(1) \right).
\]
This yields
\be
\|T_{A_\delta}\phi_\delta\|=o(1)\;\;{\rm and}\;\;\|\Xi_\delta\phi_\delta\|=o(\delta^{-2}).
\label{eq: est}
\ee
Set ${\dir}_\infty=Z_{A,\infty}+T_A$. We know have 
\[
{\dir}_\infty \phi_\delta = \psi_\delta :=\Xi_\delta\phi_\delta +T_{A_\delta}\phi_\delta
-R\phi_\delta
\]
The  estimates (\ref{eq: error}) and (\ref{eq: est}) show that
$\|\psi_\delta\|=o(1)$. The elliptic estimates applied to the {\em fixed} elliptic operator
${\dir}_\infty$ show that
\[
\|\phi_\delta\|_{1,2}\leq  C(\|\phi_\delta\| +\|\psi\|_\delta).
\]
Thus $\phi_\delta$ is bounded in the $L^{1,2}$-norm so a subsequence  (still
denoted by $\phi_\delta$) will converge strongly in $L^2$. Using again the
elliptic estimates we get
\[
\|\phi_{\delta_1} -\phi_{\delta_2}\|_{1,2}\leq  C(\|\phi_{\delta_1}
-\phi_{\delta_2}\|
+\|\psi_{\delta_1}\| +\|\psi_{\delta_2}\|).
\]
This shows the (sub)sequence $\phi_\delta$ is  Cauchy in the $L^{1,2}$ norm and
in particular it must converge in this norm to a $\phi$ which must satisfy
\[
Z_{A,\infty}\phi =T_A\phi =0
\]
due to (\ref{eq: est}).  The proof is complete. \hspace{.3cm} $\Box$

\begin{remark}{\rm  A result stronger  than stated above is true.  If the
connections $A_\delta$ are as in the adiabatic decoupling lemma and $\psi_\delta
\in L^{1,2}({\bS})$ are such that $\|\psi_\delta\|=1$ and
$\|{\dir}_\delta\psi_\delta\|=o(1)$ then a subsequence of $\psi_\delta$
converges strongly in $L^{1,2}$  to a spinor $\psi$ satisfying
(\ref{eq: adiabaticdec}).}
\label{rem: adiabatic}
\end{remark}

\bigskip

We can now conclude the proof of the theorem. From the above lemma we know that
$\phi_\delta$ converges in $L^{1,2}$.   By Sobolev inequality  the sequence $\phi_\delta$ must also converge in $L^p$, $1\leq p \leq 6$.  The last two equations in (\ref{eq: swd})  allow us to conclude that  $F_{A_\delta}$ is actually convergent in $L^3$ so   we can conclude that modulo Coulomb gauges  the connections $A_\delta$ converge {\em strongly} in  $L^{1,3}$. 

Note that $\phi$  satisfies a condition slightly weaker than (\ref{eq: 3}) namely
 \[
 \|\,|\alpha|\,\cdot \, |\beta|\, \|=0.
 \]
 We will now prove this implies (\ref{eq: 3}).  We can rewrite the limiting connection $\hat{\nabla}^A$ as a sum $\hat{\nabla}^A = \nabla^0 + B$, where  $B \in L^{1,p}({\rm End}\, ({\bS}_L))$ for  any $1\leq p < \infty$. Note that $\alpha$ and $\beta$ satisfy the elliptic Dirac equation  with  H\"{o}lder continuous coefficients
  \be
{\dir}_\infty\phi=0
\label{eq: 4}
\ee
This equation can be rewritten as
\be
{\dir}_{\nabla}^0 \phi +Q_B\phi =0
\label{eq: regular}
\ee
where $Q_B\in L^{1,p}({\rm End}\,({\bS}_L))$, $1\leq p <\infty$. Since $\phi \in L^{1,2}({\bS}_L)\cap L^\infty({\bS}_L)$ we deduce
\[
Q_B\phi \in L^{1,2}({\bS}_L).
\]
Using this in  (\ref{eq: regular}) we deduce $\phi \in L^{2,2}$.  Note that $\hat{\nabla}^{\pm}_\zeta=\nabla_\zeta^{A, \pm}\pm \lambda/2$  where $\nabla^A$ denotes the spinor connection  on ${\bS}_L$ obtained by tensoring the connection $A$ with the Levi-Civita connection defined by the metric $g_1$. Applying ${\obnaf}^+$ to the second equation in (\ref{eq: 4}) and  using the commutator equality (\ref{eq: excellent}) we
deduce after some simple manipulations that 
\[
{\obnaf} {\obnaf}^*\alpha -{\nab}^2_\zeta\alpha +|\beta|^2 \alpha
+\lambda{\obnaf}\beta=0.
\]
Since ${\obnaf}\beta =0$ and $\alpha \otimes \beta \equiv 0$ we deduce
\[
{\obnaf}{\obnaf}^*\alpha -{\nab}^2_\zeta \alpha =0.
\]
We can rewrite the above equation as
\[
{}^\flat\bar{\nabla}^0 ({}^\flat\bar{\nabla}^0)^*\alpha -{\nabla}^2_\zeta \alpha +
P(\nabla \alpha) + Q\alpha=0
\]
where the ``coefficient''  $P \in L^{1,2} \subset L^6$ while $Q\in L^2$ is at least $L^2$ (it is a linear combination of the connection $A$ and its derivatives) while the differential operator above is a generalized Laplacian. This is more than sufficient to 
apply the unique continuation principle in  Thm. 4.3 of [H] to deduce that  if $\alpha$ vanishes on an open subset of $N$  then it must vanish everywhere.

A dual argument  proves a similar result for $\beta$.  Since the product $\alpha \otimes \beta$ is identically zero we deduce that one of them must vanish on an open set and hence everywhere. The equality (\ref{eq: 3}) is proved and this completes the proof of the theorem.\hspace{.3cm}  $\Box$

\subsection{The Seiberg-Witten equations on circle bundles} 

When $N$ is the total  space of a circle bundle over a surface the above  result can be given a more precise  form.

Let $N_\ell$ be the total space of  a degree $\ell$ principal $S^1$ bundle  
\[
S^1\hookrightarrow N_\ell \stackrel{\pi}{\ra}\Sigma
\]
where $\Sigma $ is a compact surface of genus $g\geq 1$.   Fix a complex  structure on $\Sigma$ and denote by $K$ the canonical line bundle. Then
\[
{\can}\cong  \pi^* K.
\]
 Moreover, according to (\ref{eq: anisotrop2}), in this  case the limiting connection $\nabla^\infty$ on ${\can}$ coincides with the  connection pulled back from  $K$. 

Now fix a ${\rm spin}$ structure on $\Sigma$ by choosing a square root $K_\Sigma^{1/2}$ on $\Sigma$. This defines by pullback a spin structure on  $N$. Denote by $(N_\ell, \eta_\delta, h_\delta)$ the   Boothby-Wang $(K, \lambda)$  structure described  in $\S1.3$.   
 
  Fix a ${\rm spin}^c$ structure on $N$ given  by a line bundle $L_N$.   If the
   Seiberg-Witten equations $SW_{\delta}$  corresponding to  the above 
  ${\rm spin}^c$ structure have solutions for every $\delta \gg 1$  we deduce 
  from the above theorem that $L_N$ admits a connection $A$ such that 
  $F_A(\zeta ,\cdot)\equiv 0$. This implies $\pi_*F_A \equiv  0$ where $\pi_*$ denotes the integration along the fibers of $\pi:N_\ell \ra \Sigma$. From the Gysin sequence of this fibration we deduce 
  that $c_1(L_N) \in \pi^* H^2(\Sigma)$ i.e.  $L_N$ is the pullback of
   a line  bundle $L_\Sigma \ra \Sigma$.  In particular we have the following vanishing 
   result.

\begin{corollary}{\rm    Fix a  ${\rm spin}^c$ structure on $N_\ell$ which is 
not the pullback of any ${\rm spin}^c$ structure on $\Sigma$. Then  the 
Seiberg-Witten equations $SW_{\delta}$ have no solutions for 
 any $\delta \gg 1$. In particular, when $|\ell|=1$   then  for any {\em
 nontrivial} ${\rm spin}^c$ structure on $N$ and  for all $\delta \gg 1$ the equations $SW_\delta$ have no
 solutions.}
\end{corollary}

 Now fix a ${\rm spin}^c$ structure on  $N_\ell$ defined by a complex line bundle
  $L_N=\pi^*L_\Sigma$. Any of adiabatic limit of the solutions of  Seiberg-Witten  equation
 is a pair $(A, \phi)$ satisfying  the following  conditions.

 \noindent (a) A connection  $A$ on $L_N$ with horizontal curvature i.e. $F_A(\zeta, \cdot )\equiv 0$. 
 
 \noindent (b) Sections $\alpha \in {\can}^{-1/2}\otimes L_N$ and $\beta \in
 {\can}^{1/2}\otimes L_N$ such that
 \[
 Z_{A, \infty}\phi =T_A\phi=0 \Longleftrightarrow {}^\flat\bar{\nabla}^A \beta =0=\hat{\nabla}^A_\zeta \alpha\;\;{\rm and}\;\;  
 ({}^\flat\bar{\nabla}^A)^* \alpha =0 = \hat{\nabla}^A_\zeta\beta
 .
 \]
 \[
 \|\alpha\|\cdot \|\beta\|=0
 \]
where the derivation $\hat{\nabla}^A_\zeta$ is the tensor  product of the corresponding derivations on $\pi^*K^{\pm 1/2}_\Sigma$ and $L_N$.

\noindent  (c) $\frac{1}{2}(|\alpha|^2 -|\beta|^2) ={\bf
 i}F_A(\zeta_1,\zeta_2)$.  

We will refer to the conditions (a)-(c) above  as the {\em adiabatic Seiberg-Witten
equations}.  A pair $(A,\phi)$ satisfying these conditions  will be  called an {\em adiabatic solution} of the (adiabatic)  Seiberg-Witten equations.  We denote by ${\cal A}_\infty={\cal A}_\infty(\ell, g, L_N)$ the collection of gauge equivalence classes of adiabatic solutions.  Note that it is not a priori clear whether any adiabatic solution is in fact an adiabatic limit. All we can state at this point is that 
\[
{\cal S}_\infty:=\lim{\cal S}_\delta  \subset {\cal A}_\infty.
\]
 
 To make further progress  understanding the nature of the adiabatic solutions we need to distinguish two situations.
 
 \noindent {\bf A.} {\bf Reducible adiabatic solutions} i.e. $\|\alpha\|+\|\beta\| =0$. In this case  the connection $A$
 must be flat.  Thus it is uniquely defined by its holonomy representation
 \[
 {\rm hol}_A:\pi_1(N_\ell)\ra S^1.
 \]
 The fundamental group of  $N_\ell$  can be presented as
 \[
 \pi_1(N_\ell)= \lan a_1,b_1,\ldots a_g,b_g,f|
 f^{-\ell}\cdot \prod_{i=1}^g[a_i,b_i]=[a_j,f]=[b_k,f]=1,\;\;\forall j,k\ran 
 \]
 so that the space of its $S^1$ representations  is a $2g+1$-dimensional torus
 if $\ell =0$   while if $\ell \neq 0$  it is a collection of $|\ell|$ tori
 $T^{2g}$ parameterized by the group of $|\ell|$-th roots of unity. For any   representation $\rho: \pi_1(N)\ra S^1$ we have $\rho (f)=\exp (2\pi k{\bf i}/\ell)$ for some $k\in {\bZ}$  and if we denote by $L_{N,\rho}$ the line bundle it determines on $N$ we have
 \be 
 c_1(L_{N,\rho}) =  \hat{k} \in {\bZ}_{|\ell|}\subset H^2(N_\ell,{\bZ}).
 \label{eq: holonomy}
\ee
 To see this,  denote by $A_\rho$ flat the connection on  $L_{N,\rho}$ with
 holonomy $\rho$ and set $B_\rho := A_\rho + {\bf i}k/\ell\eta$. Then  the curvature of
 $B_\rho$ is purely horizontal, $F_{B_\rho}={\bf i}k/\ell d\eta$ and its holonomy
 along the fibers is zero.  Thus $B_\rho$ is the pull-back of a connection
 $D_\rho$ on  a line bundle $L_{\Sigma,\rho}\ra \Sigma$ whose  curvature  satisfies
 \[
 \pi^*F_{D_\rho} = {\bf i}k/\ell d\eta =-2{\bf i}k \ast_N \eta =-2{\bf i}k
 \pi^*d\,{\rm vol}_\Sigma.
 \]
 Thus
 \[
 c_1(D_\rho) =k/\pi \,d\,{\rm vol}_\Sigma
 \]
 i.e. $\deg L_{\Sigma, \rho} =k$ since by construction ${\rm vol}\,(\Sigma)=\pi$.   The equality (\ref{eq: holonomy}) now follows from the equality $L_{N,\rho}=\pi^*L_{\Sigma, \rho}$ combined with Gysin's exact sequence.
 To summarize, when $c_1(L_N)$ is a torsion class, the collection of gauge equivalence
 classes of reducible adiabatic solutions  forms  is a torus of dimension $2g$.
 
\bigskip
 
 \noindent {\bf B.} {\bf Irreducible adiabatic solutions} i.e. 
 $\|\alpha\|+\|\beta\|\neq 0$.
 
 \begin{lemma}{\rm If $ \|\alpha\|+ \|\beta\|\neq 0$ then (modulo gauge
 transformations) the connection $A$ is the pull back of a connection on a
 line bundle $L_\Sigma\ra \Sigma$.}
 \end{lemma}
 
 \noindent{\bf Proof of the lemma}\hspace{.2cm}   Because the connection $A$ is  
 vertically flat  one can prove easily that  the holonomy along a fiber is
 independent of the particular fiber.   Since either $\alpha$ or $\beta$ is 
 not identically zero one deduces from the  condition (b) above that this
 holonomy is trivial.  $\Box$
 
 \bigskip

  Let $(A,\alpha, \beta)$ adiabatic limits as in the above lemma. Thus there
  exists a connection $\underline{A}$ on $L_\Sigma$
  such that $\pi^* L_\Sigma= L_N$ and $\pi^*\underline{A}=A$. Since $\alpha$ and
  $\beta$ are covariant constant along fibers they can be regarded as sections
  of $L_\Sigma \otimes K_\Sigma^{\pm 1/2}$. We can decide which of $\alpha$ or $\beta$ vanishes.  More precisely if
 \[
 \deg L_\Sigma   < \deg K_\Sigma^{-1/2} = 1-g
 \]
 then the line bundle  $L_\Sigma \otimes K_\Sigma^{1/2}$ cannot admit holomorphic sections so that $\beta \equiv 0$.  In particular, $\alpha \not\equiv 0$  In this case, using the  equality (\ref{eq: chern}) we deduce after an integration along fibers that
\[
1-g >\deg L_\Sigma = \frac{\bf i}{2\pi}\int_N \eta  \wedge  F_A =\frac{1}{2}\|\alpha\|^2.
\]
This is impossible  since  $\alpha \neq 0$. If
 \[
 \deg L_\Sigma > \deg K_\Sigma^{1/2} = g-1
 \]
 then  $L_\Sigma\otimes K_\Sigma^{-1/2}$ cannot admit antiholomorphic sections so that $\alpha \equiv 0$. Reasoning as above we deduce  another contradiction so we can conclude   the adiabatic limit set is empty when $|\deg L_\Sigma | > g-1$.   

We now analyze what happens when $|\deg L_\Sigma | \leq g-1$.  We will discuss only the case $1-g \leq \deg L_\Sigma \leq 0$. The other half is completely similar.  

Using the equality (\ref{eq: chern}) again we deduce
\[
\frac{1}{4\pi}(\|\alpha\|^2-\|\beta\|^2)  =\deg L_\Sigma \leq 0.
\]
Since one of $\alpha$ or $\beta $ is identically zero  we deduce  that  $\alpha \equiv 0$. Thus the adiabatic limit set   consists of  {\em normalized holomorphic pairs} i.e pairs of the form
\[
({\rm holomorphic\; structure \; on}\; L_\Sigma\; , \; {\rm holomorphic\; section }\; \beta\; {\rm of}\; L_\Sigma\otimes K_\Sigma^{1/2})
\]
where $\beta$ is normalized by 
\be
\|\beta\|_\Sigma^2= 2 |\deg L_\Sigma|
\label{eq: normal}
\ee
 Note that when $\deg L_\Sigma = 0$ we are actually in case {\bf A}  so  it should be excluded from this discussion.  When $\deg L_\Sigma =1-g$ the bundle $L_\Sigma\otimes K_\Sigma^{1/2}$ is topologically trivial.   The only topologically trivial holomorphic line bundle  which admits a holomorphic section is the {\em holomorphically trivial} line bundle.  Hence,  when $\deg L_\Sigma =1-g$
 there exists  exactly one irreducible adiabatic solution (up to a gauge transformations). 

The above picture  resembles very much the  exact  computations in the trivial case $\ell =0$.  There are however  some notable differences.
When $\ell =0$  if the ${\rm spin}^c$ structure is non-trivial there exist no
reducible solutions for the simple fact that there exist no flat connections.
This is no longer the  case when $\ell \neq 0$. Moreover,  according to Gysin's exact sequence 
 the kernel of the morphism
\[
\pi^*: {\bZ}\cong H^2(\Sigma, {\bZ})\ra H^2(N_\ell, {\bZ})
\]
is the subgroup $\ell {\bZ}$.  If $\ell =0$ this means  the ${\rm spin}^c$-structure 
$L_N$ uniquely determines  a line bundle $L_\Sigma \ra \Sigma$  such that $\pi^*L_\Sigma \cong L_N$.

If $\ell \neq 0$  there are infinitely many line bundles  on $\Sigma$ with the above property and their degrees are congruent modulo $\ell$  with $c_1(L_N) \in {\bZ}_\ell \subset H^2(N_\ell,{\bZ})$. On the other hand,  in case {\bf B} the only line  bundles on $\Sigma$ relevant  in the adiabatic limit are those of degree in the interval $[-(g-1), (g-1)]$.

Assume now $|\ell | \geq 2g$  and 
\[
c_1(L_N) \in  {\cal V}_{g,\ell}= \{ g,g+1, \cdots ,|\ell |-g
\}  \;({\rm mod}\; \ell)\subset {\bZ}_\ell 
\]
Then  there is no integer $n$ such that $n\;{\rm mod}\; \ell \in {\cal V}_{g, \ell}$ and $0<|n|\leq g-1$.  In other words if $c_1(L_N)\in {\cal V}_{g,\ell}$ then there exists no line bundle $L_\Sigma \ra \Sigma$  of degree $0<|\deg L_\Sigma|\leq g-1$ such that $\pi^*L_\Sigma =L_N$. We have  thus proved the following result.

\begin{corollary}{\rm (a) If $|\ell| \geq 2g \geq 4$  and $L_N\ra N_\ell$ is a line 
bundle on $N_\ell$ with $c_1(L_N)\in {\cal V}_{g, \ell}$ then the collection  ${\cal A}_\infty$ of adiabatic  solutions corresponding to the ${\rm spin}^c$ structure 
defined by $L_N$ consists only of reducible configurations i.e.  pairs 
\[
({\rm flat\;\;connections,\;\;zero\;\;spinor}).
\] 
(b) If $g=1$ i.e. $\Sigma$ is a torus,    then  the set ${\cal A}_\infty$ of  adiabatic
 solutions  corresponding  to {\em any} pulled-back  ${\rm spin}^c$ structure consists only of reducible configurations.}
\label{cor: reduce}
\end{corollary}

 Any reducible adiabatic solution is obviously  a bona-fide solution  of $SW_\delta$ for all $\delta >0$.  Thus for the ${\rm spin}^c$ structures in the above corollary we have trivially
\[
{\cal S}_\infty ={\cal A}_\infty
\]
In these cases we  have ${\cal S}_\infty \subset {\cal S}_\delta$  and   it is thus natural to ask  whether  this  is a {\em strict} inclusion for {\em all} $\delta \gg 1$ i.e. whether  there  exist  adiabatic solutions which are  limits of {\em
irreducible} solutions. The estimating techniques of [GM] (see also [Ber]) show that for sufficiently large $\delta$ if $(A,\phi)$ is a
solution of the Seiberg-Witten equation $SW_\delta$  on a circle bundle over a
{\em flat torus} then $\dim \ker {\dir}_A \leq 2$. This seems to indicate that the chances  that in this case nontrivial solutions exist are very slim.   More generally  we can ask whether ${\cal S}_\infty$ does indeed coincide with ${\cal A}_\infty$ and if this is the case  how can one use the explicit knowledge of  ${\cal A}_\infty$ to obtain  better informations  about ${\cal S}_\delta$ with large $\delta$.

 A first step in this direction  is contained in the following  result.

\begin{proposition}{\rm  Fix a ${\rm spin}^c$ structure on $N_{g, \ell}$ defined by  a 
line  bundle $L_N \ra N$ such that $c_1(L_N)\in {\cal V}_{g,\ell}$.  Then} 
\[
{\cal A}_\infty(g, \ell, L_N)={\cal S}_\delta (g, \ell, L_N),\;\;\;\forall \delta \gg 1.
\]
\label{prop: adiabatic}
\end{proposition}

The range ${\cal V}_{g, \ell}$ of ${\rm spin}^c$ structures will be called the
 {\em adiabatically stable range}.
 
 \bigskip

\noindent{\bf Proof} Assume the contrary. Thus,  for all $\delta \gg 1$ there
exists $(A_\delta, \phi_\delta)\in {\cal S}_\delta$ such that $\phi_\delta \neq
0$. Set $\psi_\delta=\frac{1}{\|\phi_\delta\|}\phi_\delta$ so that
$\|\psi_\delta\|=1$, $\forall \delta$.  We can now apply the adiabatic
decoupling lemma to the sequence $(A_\delta, \psi_\delta)$ to conclude that
there exists  $\psi \in L^{1,2}({\bS}_L)$ such that  $\|\psi\|=1$ and
\[
Z_{A,\infty}\psi=T_A\psi =0
\]
where $A=\lim A_\delta$.    The discussion preceding Corollary \ref{cor:
reduce} shows that this thing is impossible when $c_1(L) \in {\cal
V}_{g,\ell}$. The proposition is proved.  \hspace{.3cm}  $\Box$

\section{Reversing the adiabatic process}

In this final section we  describe a  range of ${\rm spin}^c$ structures  on a
circle bundle over a surface for which ${\cal A}_\infty =\lim_{\delta\ra
\infty}{\cal S}_\delta$.   This boils down to constructing   irreducible
solutions of the $SW_\delta$ equations starting from irreducible adiabatic
solutions. 

\subsection{Statement of the main result}

We continue to use the same notations as in \S 3.3.  Fix a line bundle
$L_\Sigma \ra \Sigma$ and  set $L=\pi^*L_\Sigma\ra N$. 

For $|\ell|\geq 2g-1\geq 0$ define
\[
{\cal R}_{\ell, g}=\{0\}\cup \{g-1, \cdots , |\ell|-g+1\}\subset {\bZ}_\ell.
\]
The main result of this section is the following.

\begin{theorem}{\rm If  $|\ell|\geq 2g-1$ and $c_1(L)\in {\cal R}_{\ell, g}$
then
\[
\lim_{\delta \ra \infty} {\cal S}_\delta(g,\ell, L)={\cal A}(g,\ell, L)
\]
i.e. for any $(\phi, A)\in {\cal A}_\infty$ there exists a solution
$(\phi_\delta, A_\delta)$ of $SW_\delta$ which (modulo ${\gG}_L$ converges in
$L^{1,2}$-norm to $(\phi,A)$.}
\label{th: reverse}
\end{theorem}

Note that when $|\ell|\geq 2g$ and $c_1(L)\neq \pm (g-1)\;{\rm (mod)}\;\ell$ this result is clear
since in this case ${\cal A}_\infty $ consists only of reducible solutions so
that trivially
\[
{\cal A}_\infty\subset {\cal S}_\delta,\;\;\forall \delta.
\]
Thus we have to consider only the case $|\ell|\geq 2g-1$ and $c_1(L)=\pm (g-1)$
mod $\ell$.  Under this assumption the moduli space of adiabatic solutions
consists of two components. A torus of reducible solutions and a unique
irreducible solution $(A_0, \phi_0)$. The reducible solutions are obviously
adiabatic limits so we only have to prove that $(A_0, \phi_0)$ can be
approximated by solutions of $SW_\delta$, $\delta \gg 1$.

We will  consider only the case $c_1(L)=-(g-1)$.  The assumption $|\ell|\geq
2g-1$ shows the {\em unique} line bundle $L_\Sigma \ra \Sigma $ such that
$\pi^*L_\Sigma=L$ is $K^{-1/2}_\Sigma$. Thus $L={\can}^{-1/2}$. In this case
\[
{\bS}_L\cong {\can}^{-1}\otimes {\uc}
\]
This unique (modulo ${\gG}_L$) irreducible adiabatic solution  is  $(\phi_0,A_0)$ where
$A_0$ is the pullback of the Levi-Civita connection on $K_\Sigma^{-1/2}$ while
 $\phi =0 \oplus c_0$ where $c_0$ is a complex constant normalized by 
\be
\frac{1}{2}|c_0|^2 =-{\bf i}F_{A_0}(\zeta_1, \zeta_2).
\label{eq: good}
\ee
Modulo a constant gauge transformation we may assume $c_0$ is a {\em real, positive} constant.
 
 Before we begin the proof of Theorem \ref{th: reverse} we want to make a few
 elementary observations.

 For any connection $A$ on $L$ define
\[
{\dir}_{A,\infty}=Z_{A,\infty}+T_A.
\]
When $A$ is the connection $A_0$ discussed above we will write simply
${\dir}_\infty$.   The equality (\ref{eq: anisotrop1}) implies
\[
{\dir}_{A,\delta}=\delta Z_{A,\infty} +T_A +\lambda_\delta/2.
\]
In particular 
\be
{\dir}_\delta\phi_0= \lambda_\delta/2 \phi_0
\label{eq: almost}
\ee
where ${\dir}_\delta:={\dir}_{A_0,\delta}$. 

 Note also that  if $\psi =\psi_0\oplus \psi_1 \in  \ker{\dir}_\infty$ is such that  $Z_{A_0,\infty}\psi =T_{A_0}\psi =0$, $\psi_0\equiv 0$ and
\[
\int_N {\re}\lan \psi, {\bf i} \phi_0\ran\, dv_1 =0
\]
then $\psi$  must be  a {\em real} multiple of  $\phi_0$.   This simple
observation will play a   crucial role in estimating the eigenvalues of the 
linearizations of $SW_\delta$.

\subsection{The set-up}

An important part  of the proof of Theorem \ref{th: reverse} is the
construction of a suitable functional framework. This is what will be
accomplished in this subsection.

We fix $A_0$ as a reference connection. Any other connection on $L$ will have
the form $A=A_0+a$, $a\in {\bf i}\Omega^1(N)$. Set for simplicity
\[
F_0=F_{A_0},\;\;F_a= F_{A_0+a},\;\;,{\dir}_\delta={\dir}_{A_0,\delta},\;\;{\dir}^\delta_a={\dir}_{A_0+a,\delta}.
\]
The configuration space for the $SW$-equations is
\[
{\cal C}=L^{1,2}({\bS}_L\oplus {\bf i} T^*N).
\]
The elements of ${\cal C}$ will be generically denoted by ${\gc}=(\psi,a)$.
This configuration space is acted upon by the gauge group ${\gG}_L$ consisting of $L^{2,2}$-gauge
transformations of  $L$. As usual, a $\delta$-subscript will indicate  the
corresponding object is defined in terms of the deformed metric $g_\delta$. The
solutions of $SW_\delta$ are critical points of the  functional
\[
{\gf}_\delta:{\cal C}_\delta\ra {\bR},\;\;{\gf}_\delta(\psi,a)=\frac{1}{2}\int_Na \wedge (
F_0+F_a)+\frac{1}{2}\int_N\lan\psi, {\dir}_a^\delta\psi\ran dv_\delta.
\]
Lemma \ref{lemma: grad} shows that the $L^2$-gradient of ${\gf}_\delta$ is 
\[
\nabla {\gf}_\delta\!\mid_{(\psi,a)}= {\dir}_a^\delta\psi \oplus {\bf c}_\delta^{-1}(\tau_\delta(\psi))-\ast F_a).
\]

Define
\[
{\cal W}_\delta =\{ a\in L^{1,2}({\bf i}T^*N)\; ;\; d_\delta^*a=0\}
\]
and denote by $\bar{\cal W}_\delta$ its $L^2$-closure. Hodge theory shows that any orbit of ${\gG}_L$ in ${\c}$ intersects  ${\b}_\delta = L^{1,2}({\bS}_L)\oplus \w_\delta$  along a discrete set  so that  in detecting the critical points of  ${\gf}_\delta$ it suffices to  study the restriction of ${\gf}_\delta$ to ${\b}_\delta$.  Since the orbits of ${\gG}_L$ do not intersect ${\b}_\delta$ orthogonally the  $L^2$-gradient of the restriction to ${\b}_\delta$ is
\[
\nabla {\gf}_\delta\!\mid_{(\psi,a)}={\dir}_a^\delta\psi \oplus P_\delta(\tau_\delta(\psi)-\ast_\delta F_a)
\]
where $P_\delta$ denotes the $L^2_\delta$-orthogonal projection onto $\bar{\w}_\delta$. The new  configuration space ${\b}_\delta$   has a residual $S^1$-action- the constant gauge transformations. Thus if  $(\psi, a)\in {\b}_\delta$ then
\[
\exp({\bf i}\theta)(\psi, a)=\exp({\bf i} \theta)\cdot \psi, a).
\]
The fixed point set of this action is $0\oplus \w_\delta$. Set 
\[
{\b}^*_\delta=\{(\psi, a)\in {\b}_\delta\; ; \; \psi \neq 0\}.
\]
The  tangent space to the  $S^1$-orbit through ${\gc}=(\psi, a)$ is ${\cal O}_{\sgc}={\bR}({\bf i}\psi, 0)$. The $L^2$- orthogonal projection onto ${\cal O}_{\sgc}^\perp$  is 
\[
Q_{\sgc}(\dot{\psi}, \dot{a})=(\psi-\|\psi\|^{-1}_\delta{\re}\left(\int_N\lan \dot{\psi}, {\bf i}\psi \ran dv_\delta\right){\bf i} \psi \, , \,\dot{a}).
\]
We will work on the  quotient ${\b}_\delta^*/S^1$. This is a smooth  Banach manifold. The $L^2$-gradient of the functional induced by ${\gf}_\delta$ on this quotient is
\[
\nabla {\gf}_\delta\!\mid_{\sgc} =Q_{\sgc}{\dir}_a^\delta \psi \oplus P_{\delta}(\tau_\delta(\psi)-\ast_\delta F_a).
\]
To compute the linearization of this gradient we will work ``upstairs'' and then orthogonally project ``downstairs''.

On ${\c}_\delta$ the linearization of the gradient  at ${\gc}=(\psi, a)$ is
\[
{\gF}_{\sgc}(\dot{\psi},\dot{a})= ({\dir}_a^\delta +{\bf c}_\delta(\dot{a})\psi)\oplus  (\dot{\tau}_\delta(\psi, \dot{\psi})-\ast_\delta d\dot{a})
\]
where $\dot{\tau}_\delta(\psi, \cdot)$ denotes the linearization of  $\tau_\delta$ at $\psi$ and it is given by  (dropping the $\delta$ subscript for simplicity)
\[
\dot{\tau}(\psi, \dot{\psi})=\frac{1}{2}\sum_{i=0}2\left(\lan \dot{\psi},{\bf c}(\eta_i)\psi\ran  + \lan \psi , {\bf c}(\eta_i)\dot{\psi}\ran\right)\eta_i
\]
\[
={\bf i}\sum_{i=0}^2{\im}\lan\psi, {\bf c}(\eta_i)\dot{\psi}\ran \eta_i.
\]
When $\psi$ is the spinor $\phi_0 =(0,c_0)$ defined in the previous subsection we set $\dot{\tau}_\delta(\dot{\psi}):=\dot{\tau}_\delta(\phi_0, \dot{\psi})$. If we write $\dot{\psi}$ as $\dot{\psi}_0\oplus\dot{\psi}_1 \in L^{1,2}({\can}^{-1}\oplus {\uc}$ then a simple computation shows that, viewed as an endomorphism of ${\bS}_L$, $\dot{\tau}(\dot{\psi})$ is  given by
\[
\dot{\tau}(\dps)=c_0\left[
\begin{array}{cc}
-{\re}\dot{\psi}_1 &  \dps_0 \\
\dot{\bar{\psi}}_0 & {\re}\dot{\psi}_1
\end{array}
\right].
\]
The linearization of $\nabla {\gf}_\delta$ at  $[{\gc}]\in {\b}_\delta/S^1$ is 
\[
{\gL}^\delta_{\sgc}\left[ \begin{array}{c}
\dot{\psi} \\
\dot{a} 
\end{array}\right] =\left[  \begin{array}{c}
Q_{\sgc}({\dir}_a^\delta \dps  + {\bf c}_\delta(\dot{a})\psi)\\

P_\delta(\dot{\tau}_\delta(\psi, \dps) -\ast_\delta d\dot{a})
\end{array}
\right]
\]
$\forall (\dps, \dot{a})$ such that 
\[
{\re}\int_N\lan \dps, {\bf i} \psi\ran dv_g=0,\;\;d_\delta^*\dot{a}=0.
\]
Denote by ${\gL}^\delta_0$ the linearization at ${\gc}_0=(\phi_0,A_0)$ and  set  
\[
X_\delta  ={\cal O}^\perp_{{\sgc}_0}
\]
equipped with the $L_\delta^2$-norm. We can regard ${\gL}_0^\delta$ as a densely defined unbounded operator $X_\delta \ra X_\delta$. Its domain  is
\be
Y_\delta=X_\delta \cap L_\delta^{1,2}.
\label{eq: setup}
\ee
Our next result  will play a crucial role in the  proof of Theorem  \ref{th: reverse}

\begin{proposition}{\rm  (a) The unbounded operator ${\gL}_\delta $ is closed, selfadjoint and Fredholm.

\noindent (b) Set $z_\delta ={\rm dist}(0, \sigma({\gL}_\delta))$ where $\sigma({\gL}_0^\delta)$ denotes the spectrum of ${\gL}_0^\delta$. Then 
\[
 z_\delta \sim \delta^{-1}
\]
where for any two functions $f(\delta)$ and $g(\delta)$ we  set}
\[
 f\sim g \Longleftrightarrow \exists c> 1; \; c^{-1}f(\delta ) \leq g (\delta)\leq cf(\delta),\;\;\forall \delta \gg1.
\]
\label{prop: e-value}
\end{proposition}

\noindent {\bf Proof}\hspace{.3cm} The first part follows easily using the ellipticity of ${\dir}^\delta$ and of the Hodge operator $d + d_\delta^*$. The details are left to the reader.

To prove the second part we argue by contradiction.  Note first that the equality 
${\dir}_\delta \phi_0= \frac{\lambda}{2\delta}\phi_0$ implies that $z_\delta =O(\delta^{-1})$. Assume  now that for all  $\delta \gg 1$ there exists a pair $(\dps_\delta, \dot{a}\delta ) \in Y_\delta$ and $\mu_\delta\in {\bR}$ such that
\be
\left\{
\begin{array}{rcl}
Q_\delta({\dir}_\delta \dps  + {\bf c}_\delta(\dot{a})\psi) & = & \mu_\delta \dps_\delta\\

P_\delta{\bf c}_\delta^{-1}(\dot{\tau}_\delta(\dps)) -\ast_\delta d\dot{a} & = & \mu_\delta \dad
\end{array}
\right. 
\label{eq: linear1}
\ee
\be
\|\dps_\delta\| +\delta^{1/2} \|\dot{a}_\delta\|_{\delta}=1, \;\;|\mu_\delta|=z_\delta = o(\delta^{-1}).
\label{eq: normalize}
\ee
(Recall that $\|\cdot \|_\delta$ denotes the $L^2$-norm defined in terms of $g_\delta$.) Write as usual $\dps_\delta =\dps_0^\delta\oplus \dps_1^\delta $. Since 
\[
{\re}\int_N \lan \dps_\delta, {\bf i} \phi_0\ran dv_\delta =0
\]
we deduce ${\im} \dps_1^\delta =0$.   On the other hand 
\[
\int_N {\re}\lan {\bf c}_\delta(\dad)\phi_0, {\bf i}\phi_0\ran dv_\delta =0.
\]
Moreover, since ${\dir}_\delta =\delta Z_{A_0,\infty}+T_{A_0}+\lambda_\delta/2$ we deduce  that if $\dps \perp {\bf i}\phi_0$ then
\[
\int_N{\re}\lan {\dir}_\delta\dps, \phi_0 \ran dv_\delta =\frac{\lambda_\delta}{2}\int_N{\re}\lan \dps, {\bf i}\phi_0 \ran dv_\delta =0
\]
The above two equalities show that we can drop  the $Q_\delta$ term from the first equation in (\ref{eq: linear1}).

\begin{lemma}{\rm $\dps_\delta \neq 0$ for all $\delta \gg 1$.}
\end{lemma}

\noindent {\bf Proof} Assume the contrary i.e.  $\dps_\delta \equiv 0$ for a subsequence $\delta \ra \infty$. Locally we can write
\[
\dad = {\bf i} u_\delta\eta_\delta + v_\delta{\ve}-\bar{v}_\delta\bar{\ve}
\]
where $u_\delta$ is real  while $v_\delta$ is complex valued. We have
\be
{\bf c}_\delta(\dad)\phi_0=c_0\left[
\begin{array}{c}\bar{v}_\delta \bar{\ve}\\
u_\delta
\end{array}
\right].
\label{eq: cliff1}
\ee
  The  the first equation of the system (\ref{eq: linear1}) becomes
\[
{\bf c}_\delta(\dad)\phi_0=0
\]
which implies $\dad \equiv 0$. This contradicts (\ref{eq: normalize}).  $\Box$

\bigskip

The  decisive moment in the proof of  Proposition  \ref{prop: e-value} is contained in the following result.

\begin{lemma}{\rm Set  $\Psi_\delta =\|\dps_\delta\|^{-1}\dps_\delta$. Then there exists $t\in {\bR}\setminus\{0\}$  and a subsequence $\delta \ra \infty$ such that $\Psi_\delta \ra t\phi_0$ in $L^{1,2}$.}
\label{lemma: blowup}
\end{lemma}

The proof of Lemma \ref{lemma: blowup} relies on  the following auxiliary result which  is an  immediate consequence of  Proposition 3.45 of  [BGV].

\begin{lemma}{\rm For any $\dot{a}\in \w_\delta$
\[
{\dir}_\delta {\bf c}_\delta(\dot{a}) +{\bf c}_\delta(\dot{a}) {\dir}_\delta = -2\nabla^\delta_{\dot{a}}+{\bf c}_\delta (\ast_\delta d\dot{a})
\]
where $\nabla^\delta$ denotes the ${\rm spin}^c$ connection defined by $A_0$ on ${\bS}_L$ using the metric $g_\delta$ while $\nabla^\delta_{\dot{a}}$ denotes the covariant derivative along the vector field $g_\delta$-dual to $\dot{a}$. (Note that $\dot{a}$ is by definition a purely imaginary form and so will be its dual vector field.)}
\label{lemma: dircom}
\end{lemma}

\noindent {\bf Proof of Lemma \ref{lemma: blowup}} Consider the complex valued 1-form  $\omega_\delta = {\bf i}{\re}\dps_1^\delta \eta_\delta +\dbps_0^\delta -\dps_0^\delta$. Then $c_0{\bf c}_\delta(\omega_\delta)=\dot{\tau}_\delta(\dps_\delta)$ so that
\[
\int_N \lan {\bf c}_\delta P_\delta{\bf c}_\delta^{-1}\left(\dtd(\dps_\delta)\right)\, \phi_0, \dps_\delta\ran dv_\delta
\]
\[
=\int_N \lan P_\delta {\bf c}_\delta^{-1}(\dot{\tau}_\delta(\dps_\delta))\, ,\, \omega_\delta \ran dv_\delta=c_0\int_N \lan P_\delta \omega_\delta ,\omega_\delta  \ran  dv_\delta \geq 0.
\]
Using the second equation in  (\ref{eq: linear1})  
\be
{\re}\int_N\lan {\bf c}_\delta(\ast_\delta d\dad)\phi_0, \dps_\delta\ran dv_\delta \geq -{\re}\mu_\delta \int_N\lan {\bf c}_\delta(\dad)\phi_0, \dps_\delta \ran dv_\delta.
\label{eq: pozitiv}
\ee
Applying ${\dir}_\delta $  to the first equation in (\ref{eq: linear1}) equation (remember $Q_\delta$ need not be included) we deduce using Lemma \ref{lemma: dircom}
\be
{\dir}^2_\delta\dps_\delta  -{\bf c}_\delta(\dad){\dir}_\delta \phi_0 -2\nabla^\delta_{\dad }\phi_0 +{\bf c}(\ast_\delta d \dad)\phi_0  =\mu_\delta {\dir}_\delta \dps_\delta.  
\label{eq: linear2}
\ee
This equality can be further simplified using the  equalities
\[
{\dir}_\delta \phi_0 =\frac{\lambda}{2\delta}\phi_0
\]
and  since  $\nabla^\delta$ is the tensor of the spin connection on ${\can}^{-1/2}\oplus {\can}^{1/2}$ with the pullback connection $A_0$ on ${\can}^{-1/2}=\pi^*K_\Sigma^{-1/2}$ we deduce from the computations in  \S 2.4 that
\be
\nabla_{\dad}\phi_0 = c_0 \frac{\lambda}{2\delta}((1+{\bf i})\bar{v}_\delta \oplus {\bf i} u_\delta)
\label{eq: small}
\ee
Taking the ({\em real}) inner product of (\ref{eq: linear2}) with $\dps_\delta$ and integrating  by parts we deduce
\[
\|{\dir}_\delta\dps_\delta\|_\delta^2 -2{\re} \lan \nabla_{\dad}\phi_0 , \dps_\delta \ran _\delta +{\re} \lan {\bf c}_\delta(\ast d\dad)\phi_0, \dps_\delta\ran_\delta -\frac{\lambda}{2\delta}{\re} \lan{\bf c}_\delta(\dad)\phi_0, \dps_\delta \ran_\delta 
\]
\[
=\mu_\delta \lan {\dir}_\delta \dps_\delta, \dps_\delta\ran_\delta  
\]
Using the inequality (\ref{eq: pozitiv}) and the identity (\ref{eq:  small}) we deduce the following inequality ($c$ will generically denote a constant independent of $\delta$).
\[
\|{\dir}_\delta\dps_\delta\|_\delta^2-\frac{c}{\delta}\|\dad\|_\delta \cdot \|\dps_\delta\|_\delta -|\mu_\delta|\cdot  \|\dad\|_\delta \cdot \| \dps_\delta\|_\delta \leq  |\mu_\delta| \cdot \|{\dir}_\delta \dps_\delta\|_\delta  \cdot \|\dps_\delta\|_\delta.
\]
This inequality implies
\[
\|{\dir}_\delta \dps_\delta\|_\delta =O(\delta^{-1}\|\dps\|_\delta)
\]
which is equivalent to
\[
\|{\dir}_\delta \Psi_\delta\|= O(\delta^{-1}),\;\; (\|\Psi_\delta\|=1).
\]
Using  the adiabatic decoupling lemma (or rather the remark following its proof) we deduce that on a subsequence $\Psi_\delta$ converges in $L^{1,2}$ to a {\em nontrivial} element $\Psi_\infty \in \ker Z_{A_0,\infty}\cap \ker T_{A_0}$ such that ${\re}\lan \Psi_\infty, {\bf i}\phi_0\ran =0$. Lemma \ref{lemma: blowup} now follows using the observation at the end of \S 4.1.  $\Box$

\bigskip

We can now  complete the proof of the Proposition \ref{prop: e-value}. Take the inner product with $\phi_0$ of  the first equation in (\ref{eq: linear1}) to deduce
\[
\left(\frac{\lambda}{2\delta}-\mu_\delta\right) \int_N \lan \dps_\delta, \phi_0\ran dv_1 =- \int_N \lan {\bf c}_\delta(\dad)\phi_0, \dps_\delta \ran dv_1 
\]
so that
\be
\left(\frac{\lambda}{2\delta}-\mu_\delta\right) \int_N \lan \Psi_\delta, \phi_0 \ran dv_1 =- \int_N\lan {\bf c}(\dad)\phi_0, \Psi_\delta\ran dv_1.
\label{eq: est1}
\ee
By taking the inner product  with $\dad \in \w_\delta$ of  the second equation in  (\ref{eq: linear1})  we get after integrating by parts  and using (\ref{eq: normalize})
\[
\left|\int_N\lan P_\delta {\bf c}_\delta^{-1}(\dtd(\dps_\delta)),\dad\ran_\delta dv_1\right|=\left|\int_N\lan {\bf c}_\delta^{-1}(\dtd(\dps_\delta)), \dad\ran_\delta dv_1\right| =O(|\mu_\delta|  \delta \|\dad\|^2_\delta)=O(|\mu_\delta|).
\]
A simple computation shows that
\[
\left|\int_N\lan {\bf c}_\delta^{-1}(\dtd(\dps_\delta)), \dad \ran_\delta dv_1\right| =\left|\int_N\lan {\bf c}(\dad)\phi_0,
\dps_\delta\ran dv_1\right|
\]
Hence 
\[
\left|\int_N\lan {\bf c}(\dad)\phi_0, \Psi_\delta)dv_1\right|= O(|\mu_\delta|)=o(\delta^{-1})
\]
Using  this last estimate  in (\ref{eq: est1}) we deduce
\[
\left(\frac{\lambda}{2}-\delta \mu_\delta\right) \int_N \lan \Psi_\delta, \phi_0 \ran dv_1 =o(1)
\]
We  now let $\delta\ra \infty$ in the above equality.  Lemma  \ref{lemma: blowup}  leads to a contradiction. The proposition is proved.  $\Box$

\subsection{Proof of the main result}

As was indicated   at the beginning of  \S 4.1  we have to show that for all  $\delta \gg  1$ there exists a solution $(\phi_\delta, A_\delta)$ of $SW_\delta$ which converges (modulo ${\gG}_L$) to  $(\phi_0, A_0)$ as $\delta \ra \infty$. To produce such solutions we will use the technique pioneered by Taubes   in [T]  where he  proved the existence of  self-dual  connections. The abstract  result  behind this technique is contained in the following   version of the inverse function theorem. We will state it in the simplest context of maps between two Banach spaces. The extension to Banach manifolds and vector bundles  is only notationally more complicated.

\begin{lemma}{\rm   Consider  a map  between two Banach spaces  $F:Y\ra X$  which can be decomposed as $F(y)=Ly +N(y)$ where $L$ and $N$ satisfy the following conditions.

\noindent (i) $L$ is  a linear isomorphism of Banach spaces. Set 
\[
\mu=\inf_{y\neq 0}\frac{\|Ly\|}{\|y\|}
\]
\noindent (ii) There exists $\kappa>0$ such  that 
\[
\|N(y_1)-N(y_2)\|\leq \kappa (\|y_1\|+\|y_2\|)(\|y_1-y_2\|),\;\;\forall y_1,y_2 \in Y.
\]
Set ${\ve}_0 =\|N(0)\|$ and $q=\kappa \mu^{-1}$.  Suppose that
\be
{\ve}_0 <\frac{1}{4q} < \frac{1}{2}.
\label{eq: balance}
\ee
and denote by $r=r(q, {\ve}_0)$ the smallest root of the quadratic equation
\be
qr^2-r+{\ve}_0=0.
\label{eq: distance}
\ee
Then there exists  an unique $y_0\in Y$ such that}
\[
\|y_0\|\leq r(q, {\ve}_0)\;\;{\rm and}\;\; F(y_0) =0.
\]
\label{lemma: implicit}
\end{lemma}

The idea of this lemma  is intuitively clear. If $0$   is an ``almost'' solution  of the equation $F(y)=0$  (${\ve}_0$ quantifies the attribute ``almost'') and the linearization $L$ is invertible then  we can perturb $0$ to a genuine solution provided the  norm  of the inverse is sufficiently small (the attribute ``sufficiently small'' is clarified in the inequality  (\ref{eq: small}) which also incorporates an interaction with the nonlinear term). 

The proof  is an immediate application of the Banach fixed point theorem  and  can be safely left  to the reader.

In our case  the equation we want to solve is 
\[
\nabla {\gf}_\delta =0
\]
near ${\gc}_0$ and $\nabla{\gf}_\delta$ is rather a section of an infinite dimensional vector bundle.  Trivializing (via parallel transport for example) near this configuration we can regard it as an equation of the form
\[
{\cal F}_\delta (y) =0,\;\;\;y\in Y_\delta
\]
where $X_\delta$ and $Y_\delta$ are defined in (\ref{eq: setup}) and ${\cal F}: Y_\delta \ra X_\delta$.  The role  of the operator $L$ in Lemma \ref{lemma: implicit} is played by  the linearization ${\gL}_0^\delta$ this time viewed as a {\em bounded} operator $Y_\delta \ra X_\delta$.   As usual $c$ will generically denote a constant independent of $\delta \gg 1$. The first thing we want to prove is the following.

\begin{lemma}{\rm There exists $c>0$  such that}
\[
\|{\gL}_0^\delta(\psi, a)\|_{2,\delta} \geq \frac{c}{\delta}(\|\psi\|_{1,2,\delta}+\|a\|_{1,2,\delta}),\;\;\forall (\psi,a)\in Y_\delta.
\]
\label{lemma: e-value}
\end{lemma}

\noindent {\bf Proof} \hspace{.3cm} Note that  there exists $c>0$ such that
\[
\|{\gL}_0^\delta (\psi, a)\|^2_{2,\delta} +\|\psi\|^2_{2,\delta} +\|a\|^2_{2,\delta} \geq c\left(\|{\dir}_\delta \psi\|^2_{2,\delta} + \|da\|^2_{2,\delta} \right).
\]
The Weitzenb\"{o}ck formula for  ${\dir}_\delta$ coupled with the boundedness of $F_{A_0}$ and of the scalar curvature of $N_\delta$ show that
\[
\|\psi\|_{1,2,\delta} \leq c\left( \|{\dir}_\delta \psi\|_{2,\delta} +\|\psi\|_{2,\delta}\right).
\]
Similarly, the Bochner formula for the Hodge operator on 1-forms coupled with the  boundedness of the Ricci curvature shows  that there exists $c>0$ such that
\[
\|a\|_{1,2,\delta} \leq c\left( \|da\|_{2,\delta}+\|a\|_{2,\delta}\right),\;\;\;\forall a\in \w_\delta.
\]
Putting all the above together we deduce there exists $c>0$ such that
\[
\|\psi\|_{1,2,\delta} +\|a\|_{1,2,\delta} \leq c\left( \|{\gL}_0^\delta(\psi,a)\|_{2,\delta} +\|a\|_{2,\delta} + \|\psi\|_{2,\delta}\|\right).
 \]
The lemma now follows from Proposition \ref{prop: e-value}.   $\Box$

\bigskip

The  nonlinear  part  ${\cal N}_\delta$ in the Seiberg-Witten equation comes only  from the quadratic  term $\tau_\delta$.    Thus in this case we have (via  the H\"{o}lder  inequalities)
\[
\|{\cal N}(\psi_1, a_1)-{\cal N}_\delta(\psi_2, a_2)\|_{2,\delta} \leq  c(\|\psi_1\|_{4,\delta}+\|\psi_2\|_{4,\delta})\cdot \|\psi_1-\psi_2\|_{4,\delta}.
\]
At this point  the controlled  geometry of the  deformation  $g \ra g_\delta$ plays a magical role.  More precisely, since the Ricci curvature is bounded, ${\rm vol}\,(N_\delta)\sim \delta^{-1}$ and ${\rm diam}\,(N_\delta) \sim 1$ we deduce from the Sobolev estimates of  [Be], Appendix VI  that there exists $c>0$ such that for  any $f\in L^{1,2}_\delta(N_\delta)$ we have the {\em sharp} Sobolev inequality
\[
\|f\|_{4,\delta}\leq c\delta^{-1/4}(\| d f\|_{2,\delta}+\|f\|_{2,\delta})
\]
Kato's  inequality now implies immediately that for every $\psi \in L^{1,2}_\delta({\bS}_L)$ we have
\[
\|\psi\|_{4,\delta}\leq c\delta^{-1/4}\|\psi\|_{1,2,\delta}.
\]
Putting all the above together we deduce
\[
\|{\cal N}(\psi_1, a_1)-{\cal N}_\delta(\psi_2, a_2)\|_{2,\delta} \leq  c \delta^{-1/2}(\|\psi_1\|_{1,2,\delta}+\|\psi_2\|_{1,2,\delta})(\|\psi_1 -\psi_2\|_{1,2,\delta}).
\]

 Finally,   using (\ref{eq: good}) and (\ref{eq: almost}) we deduce $\|\nabla {\gf}_\delta({\gc}_0)\|_{2,\delta}\sim \delta^{-3/2}$.

We now want  to apply  Lemma \ref{lemma: implicit} with 
\[
{\ve}_0\sim \delta^{-3/2}, \;\;\mu=O(\delta^{-1})\;\;{\rm  and}\;\; \kappa =O(\delta^{-1/2}).
\]
The only problem is that the solution it postulates might be a reducible one.  We need to eliminate this possibility.  This is where  sharp asymptotics of the solutions of (\ref{eq: distance}) are needed. The quadratic formula shows
\[
r(q,{\ve}_0)= \frac{1 -\sqrt{1-4q{\ve}_0}}{4q}
\]
where ${\ve}_0(\delta)\sim \delta^{-3/2}$ . Thus
\[
r(q,{\ve}_0)\sim {\ve}_0(\delta)\sim \delta^{-3/2}.
\]
On the other hand the distance from $\phi_0$ to the reducible set  $0\oplus \w_\delta$ is $\|\phi_0\|_{2,\delta}\sim \delta^{-1/2}$.  Hence for $\delta \gg 1$ the solution detected with the help of  Lemma \ref{lemma: implicit} is  too close to $(\phi_0,A_0)$ to be reducible. Theorem \ref{th: reverse} is proved.  $\Box$


\end{document}